\documentclass[a4paper]{article}
\usepackage[left=4cm, right=4cm, bottom=3.5cm, top=3.5cm]{geometry}

\usepackage{graphicx}
\usepackage{amsmath}
\usepackage{amsfonts}
\usepackage{amssymb}
\graphicspath{{figures/}}

\newtheorem{definition}{Definition}
\providecommand{\keywords}[1]{\textbf{Keywords: } #1}

\newcommand{\Rbb}{\mathbb{R}}
\newcommand{\Sbb}{\mathbb{S}}
\newcommand{\Tbb}{\mathbb{T}}

\title{The phase space geometry underlying roaming reaction dynamics
}

\author{Vladim{\'i}r Kraj{\v{n}}{\'a}k\footnote{School of Mathematics, University of Bristol, University Walk, Bristol BS8 1TW, UK \& Johann Bernoulli Institute, University of Groningen, Nijenborgh 9, 9747 AG Groningen, The Netherlands}~~and~Holger Waalkens\footnote{Johann Bernoulli Institute, University of Groningen, Nijenborgh 9, 9747 AG Groningen, The Netherlands}}

\date{}

\begin{document}

 \maketitle

 \begin{abstract}
  Recent studies have found an unusual way of dissociation in formaldehyde. It can be characterized by a hydrogen atom that separates from the molecule, but instead of dissociating immediately it roams around the molecule for a considerable amount of time and extracts another hydrogen atom from the molecule prior to dissociation. This phenomenon has been coined  roaming and has since been reported in the dissociation of a number of other molecules. In this paper we investigate roaming in Chesnavich's CH$_4^+$ model. During dissociation the free hydrogen must pass through three phase space bottleneck for the classical motion, that can be shown to exist due to unstable periodic orbits. None of these orbits is associated with saddle points of the potential energy surface and hence related to transition states in the usual sense. We explain how the intricate  phase space geometry influences the shape and intersections of invariant manifolds that form separatrices, and establish the impact of these phase space structures on residence times and rotation numbers. Ultimately we use this knowledge to attribute the roaming phenomenon to particular heteroclinic intersections.

  \keywords{reaction dynamics, roaming, transition state theory}
 \end{abstract}

 \section{Introduction}\label{sec:introch4}
 \subsection{Roaming}
  For a long time it was believed that dissociation of molecules can only happen in two ways. Firstly, the original molecule can dissociate into smaller molecules and this is sometimes referred to as dissociation via the molecular channel. In order to dissociate, the system has to pass over a potential barrier representing the energy needed to break existing bonds and form new ones. Quantitative results on dissociation rates (or reaction rates in general) can be obtained via transition state theory. Alternatively an individual atom, called free radical, can escape from a molecule without forming new bonds and thus without passing over a potential barrier \cite{Bowman2014}. This is sometimes referred to as dissociation via the radical channel. Dissociation via both channels is well understood.
  
  Recently, however, van Zee et al. \cite{vanZee1993} reported having experimentally observed dissociation of formaldehyde (H$_2$CO) with CO in low rotational levels at an energy where dissociation through the molecular channel should have rather resulted in high rotational states of CO. The two proposed explanations for this behaviour are that either at least one of the vibrational modes of the transition state is quite anharmonic or there have to be two distinct molecular channels.
  
  Townsend et al. \cite{Townsend2004} discovered in their study of formaldehyde (H$_2$CO) a new form of dissociation that appears not to be associated with the molecular or the radical channel. In the process, an H atom separates from the molecule following the radical channel, but instead of dissociating it spends a considerable amount of time near the molecule and eventually abstracting the remaining H atom from the molecule to form H$_2$. This type of dissociation is called \emph{roaming} due to the nature of behaviour of the escaping H atom. No potential barrier or dynamical transition state is known to be involved in roaming.
  
  The discovery of roaming stimulated extensive studies of formaldehyde photodissociation and roaming has since been accepted as the cause for the phenomenon observed by van Zee et al. \cite{vanZee1993}.  
  
  Bowman and Shelper \cite{Bowman2011} have studied the dynamics of H$_2$CO and CH$_3$CHO to find evidence that roaming is more connected to the radical rather than the molecular channel. At the same time, roaming was observed at energies below the radical threshold.

      
 \subsection{Known results} 

  In recent years dynamical systems theory has made an impact in chemistry by providing the means to understand the classical phase space structures that underlie reaction type dynamics \cite{Wigginsetal01,Uzeretal02,Waalkens04}. This concerns the phase space geometry that governs transport across a saddle equilibrium point referred to as the molecular channel above. These ideas make precise the notion of a \emph{transition state} which forms the basis of computing reaction rates from \emph{Transition State Theory} which can be considered to be the most important approach to compute reaction rates in chemistry \cite{Wigner37}.
  
  It is shown that surfaces of constant energy contain for energies above the saddle an invariant manifold with the topology of a sphere. This sphere is unstable. More precisely it is a normally hyperbolic invariant manifold (NHIM) which is of codimension 2 in the energy surface. The NHIM can be identified with the transition state: it forms an unstable invariant subsystem located between reactants and products. What is crucial for the computation of reaction rates is that the NHIM is spanned by the hemispheres of another higher dimensional sphere which is of codimension 1 in the energy surface and which is referred to as a \emph{dividing surface}. It divides the energy surface into a reactants region and  a products region in such a way that trajectories extending from reactants to products have exactly one intersection with one of the hemispheres and trajectories extending from products to reactants have exactly one intersection with the other hemisphere. The construction of such a \emph{recrossing-free} dividing surface is crucial for Transition State Theory where reaction rates are computed from the flux through a dividing surface which is computationally much cheaper than sampling trajectories.
  
  For the global dynamics, the NHIM is significant because of its stable and unstable manifolds. The latter have the topology of spherical cylinders or `tubes' which are of codimension 1 in the energy surface and hence have sufficient dimensionality to act as separatrices. In fact the stable and unstable manifolds separate the reactive trajectories from the non-reactive ones. The geometry of the stable and unstable manifolds and their location and  intersections in the reactants and products regions carry the full information about the transition process including, e.g., state specific reactivity \cite{EzraWaalkensWiggins2009}. In the case of two degrees of freedom the NHIM is the Lyapunov periodic orbit associated with the saddle equilibrium point and the approach reduces the\emph{ periodic orbit dividing surface} (PODS) introduced earlier by Pechukas and  Pollak and others \cite{PechukasMcLafferty73,PechukasPollak78}.

  Explaining the roaming phenomenon poses a new challenge to dynamical systems theory. The first attempts to use methods  from dynamical systems theory related to ones underlying transition state theory to explain roaming can be found in the work of Maugui\`ere et al. \cite{Mauguiere2014} in which they identified the region of the classical  phase space where roaming occurs with the aid of numerous invariant phase space structures for Chesnavich's CH$_4^+$ dissociation model \cite{Chesnavich1986}. They also introduced a classification of trajectories present in the system and matched them with the experimentally observed behaviour. The definition of roaming was formulated using the number of turning points in the radial direction in the roaming region. In light of a gap time analysis of  Maugui\`ere et al. in \cite{Mauguiere2014b}, the definition was refined by means of the number of crossings of a dividing surface constructed in the roaming region. This refined definition is the best dynamical description of roaming to date.
  
  In \cite{Mauguiere2015}, Maugui\`ere et al. studied the model of formaldehyde to find unstable periodic orbits in the roaming region. The homoclinic tangle of one such orbit was shown to be responsible for transport between two potential wells in a process that is closely linked to roaming. The  periodic orbit involved do not arise as the Lyapunov periodic orbits associated with a saddle equilibrium points and the situation is hence different from the usual setting of transition state theory  built on a potential saddle. Yet a recrossing free dividing surface can be constructed from such a periodic orbit. Such dividing surfaces may be other than spherical. It was shown in \cite{MacKay2014} and \cite{MacKay2015} that a spherical dividing surface near an index-$1$ critical point may bifurcate into a torus. In \cite{Mauguiere2013} a toric dividing surface was constructed near an index-$2$ critical point.
  
  The authors of \cite{Mauguiere2016} found that the local geometry of the energy surface in an O$_3$ model may be toric and constructed a toric dividing surface using two unstable periodic orbits.
  
  Even more recently, Huston et al. \cite{Huston2016} report that they have found a correlation between the distribution of internal energies of CO and H$_2$ with the molecular channel and roaming. Particularly at lower energies, roaming trajectories have significantly more energy in H$_2$. With increasing energy the differences decrease. Their definition of roaming is slightly different though, it involves rotation of H$_2$ around CO at a `slowly varying and elongated distance'. The precise definition involves technical conditions on H$_2$ vibrational energy, time spent at a certain minimal distance from CO with low kinetic energy and large H-H bond length.


 \subsection{Objectives and outline of this paper}
    
  Because there is no single generally accepted definition of roaming, there is a clear need for a deeper understanding of the mechanisms behind dissociations.
  
  In this work we present a detailed study of dissociation in the CH$_4^+$ model by \cite{Chesnavich1986}. We discuss all types of dynamics present in this model and explain their connection to the underlying phase space geometry and invariant structures. We construct various surfaces of section and from the dynamics on these surfaces we deduce the role of invariant manifolds in slow dissociation and ultimately show a certain structure of heteroclinic tangles that causes roaming.
  
  From the point of view of transition state theory we address two interesting problems. 
  Firstly, it is not very well understood what happens in case reactants and products are divided by multiple transition states in series, which is a problem we address in this work. Secondly, we study the role of the local energy surface geometry in interactions of multiple transition states.
  
  In the study we employ surfaces of section, all of which satisfy the Birkhoff condition \cite{Birkhoff27} of being bounded by invariant manifolds. Using the surfaces of section we can observe dynamical behaviour such as roaming, but to understand the role of the local energy surface geometry and its implications to roaming, we generalise the Conley-McGehee representation \cite{Conley1968,McGehee1969,MacKay1990} and study the dynamics on the energy surface in full $3$ dimensions.

  The paper is organized as follows. In Sec.~\ref{sec:setup} we introduce the Chesnavich's CH$_4^+$ model. 
  In Sec.~\ref{sec:dynamics} we discuss various periodic orbits and their role in setting up the problem of transport of phase space volumes between different phase space regions.
  In Sec.~\ref{sec:observations} we study the dynamics of the Chesnavich model by looking at  trajectories from various perspectives. 
  This section is followed by relating the dynamics to roaming in Sec.~\ref{sec:discussion}. 
  The invariant manifolds that govern the dynamics and in particular roaming are discussed from a global perspective in Sec.~\ref{sec:representation}.
  Conclusions are given in Sec.~\ref{sec:conclusions}


 \section{Set-up}\label{sec:setup}

 \subsection{Chesnavich's CH$_4^+$ dissociation model}\label{subsec:system}

  Like \cite{Mauguiere2014} and \cite{Mauguiere2014b}, we use the model for the $ \text{CH}_4^+ \rightarrow \text{CH}_3^+ + \text{H}$ dissociation introduced by Chesnavich \cite{Chesnavich1986}. The system is a $2$-degree-of-freedom ``phenomenological model'' that is intended for the study of multiple transition states. In this model, only one H atom is free and the CH$_3^+$ molecule is considered to be a rigid complex.

  It is a planar system and we study it in a centre of mass frame in polar coordinates $(r,\theta_1, \theta_2, p_r , p_1, p_2)$, where $r$ is the distance of the free H atom to the centre of mass (magnitude of the Jacobi vector), $\theta_1$ is the angle between a fixed axis through the centre of mass and the Jacobi vector between CH$_3^+$ and H, and $\theta_2$ is the angle representing the orientation of CH$_3^+$ with respect to the fixed axis.

  In these coordinates the kinetic energy has the form 
  $$ T = \frac{1}{2 m} \left( p_r^2 + \frac{1}{r^2} p_2^2 \right) + \frac{1}{2 I} p_1^2, $$
  where $m$ is the reduced mass of the system, and $I$ is the moment of inertia of the rigid body CH$_3^+$. The system has a rotational $SO(2)$ symmetry, which can be reduced giving a family of systems parametrised by the (conserved) angular momentum.

  This reduction can be obtained from the following canonical transformation:
  \[ \theta_1 = \theta + \psi, \quad \theta_2 = \psi, \quad p_1 = p_\theta, \quad p_2 = p_\psi - p_\theta.\]
  Then $p_\psi = p_1 + p_2=:\lambda$ is the total angular momentum and it is conserved. It follows that
  \begin{align*}
  H(r,\theta,p_r, p_\theta;\lambda) &= \frac{1}{2 m} p_r^2 + \frac{1}{2 I} p_\theta^2  + \frac{1}{2 m r^2} (p_\theta - \lambda)^2  + U(r,\theta) \\
  &= \frac{1}{2 m} p_r^2 + \frac{1}{2} \left(\frac{1}{I} + \frac{1}{m r^2} \right) p_\theta^2  - \frac{\lambda}{m r^2} p_\theta + \frac{\lambda^2}{2 m r^2} + U(r,\theta),
  \end{align*}
  where $U(r,\theta)$ is the potential energy from \cite{Chesnavich1986} that we will discuss later in Section \ref{subsec:ch4potential}. In the last expression the term $\frac{\lambda}{m r^2} p_\theta$ gives rise to a Coriolis force in the equations of motion.


 \subsection{General setting}\label{subsec:general setting}
  As explained by \cite{Bowman2014}, systems exhibiting roaming have a potential well for a small radius, representing the stable molecule, and with increasing distance between the dissociated components converges monotonously to a certain base energy, which we can assume to be $0$. This is unlike the traditional bimolecular reactions that involve flux over a potential saddle.  
  As shown in \cite{Mauguiere2016}, under certain conditions the potential $U$ admits two unstable periodic orbits that are not associated with any potential which, however, form the transition state to dissociation.  We will find these orbits and use them to construct a toric dividing surface from them.  

  The argument for the existence of the periodic orbits is as follows. The dependence of the potential $U(r,\theta)$ on $\theta$ is due to the interaction between the anisotropic rigid molecule and the free atom. When $r$ is sufficiently large, the potential $U(r,\theta)$ is essentially independent of $\theta$. This is because for sufficiently large distances, the orientation of the CH$_3^+$ molecule does essentially not influence the interaction with the free atom. Let us therefore assume for a moment that $r$ is sufficiently large, so that $U$ is rotationally symmetric and we can drop $\theta$ in the argument of  $U$. The system reduced by the rotational symmetry  then has the effective potential
  $$V_{red}(r;\lambda)=\frac{(p_\theta-\lambda)^2}{2 m r^2} + U(r),$$
  where $p_\theta$ becomes a constant of motion. The reduced system admits a relative equilibrium, provided $U(r)$ is monotonous, $U(r)<0$ and $U\in o(r^{-2})$ as $r\rightarrow\infty$. Potentials of most chemical reactions, including $\text{CH}_4^+ \rightarrow \text{CH}_3^+ + \text{H}$, meet this condition. The relative equilibrium is given by $r=r_{p_\theta}$, $p_r=0$, where $r_{p_\theta}$ is the solution of
  $$ \dot {p}_r = -\frac{\partial H}{\partial r} = \frac{1}{m r^3} (p_\theta - \lambda )^2 - \frac{\mathrm{d} U}{\mathrm{d} r} = 0. $$

\begin{figure}
  \centering
  \includegraphics[width=1\textwidth]{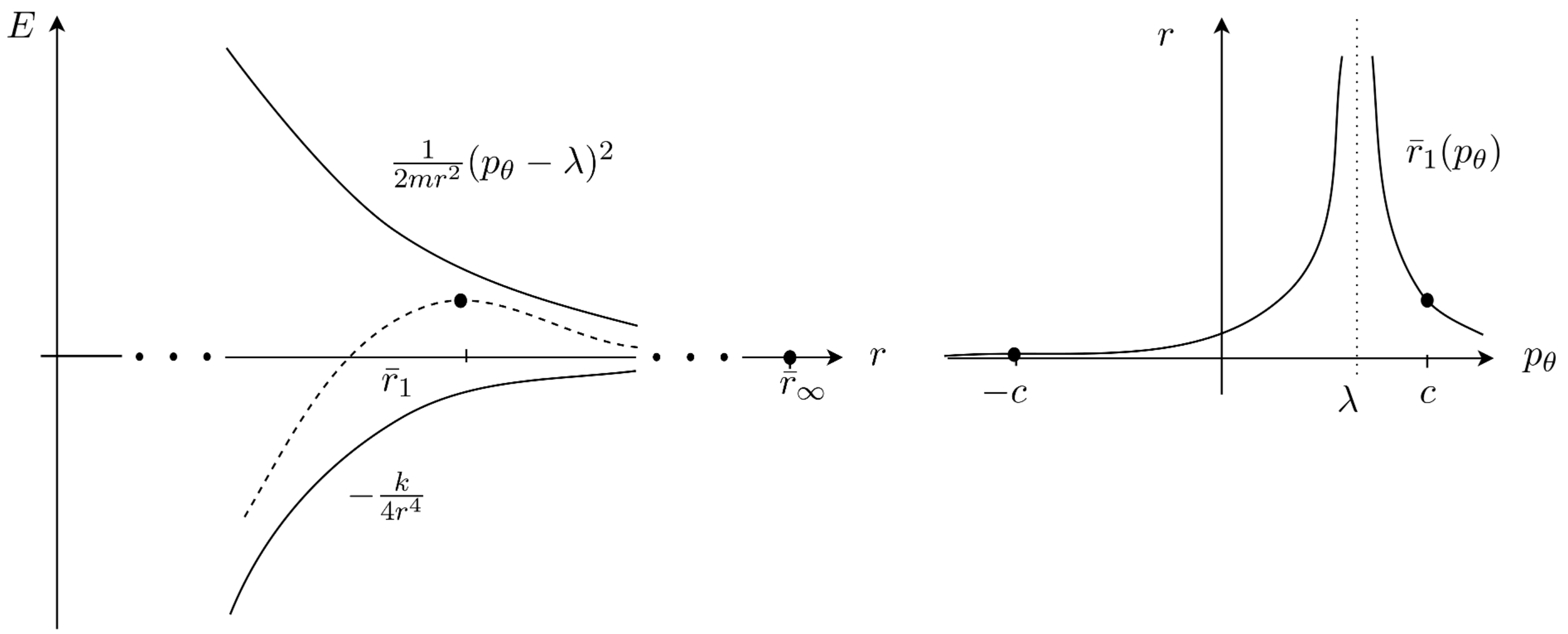}
  \caption{Schematic representation of the dominant long range potential and the ``centrifugal term'' over $r$  (left), and of $r_1$ over $p_\theta$ (right).}
  \label{fig:centrifugalTS}
  \end{figure}

  For the class of potentials $U= -c r^{-(2+\epsilon)}$ with $c,\epsilon>0$, the relative equilibrium is unstable (Fig. \ref{fig:centrifugalTS}). This follows from the reduced $1$-degree-of-freedom Hamiltonian having a saddle at this equilibrium as can be seen from computing the Hessian matrix which is diagonal with the elements
  $$\frac{\partial^2 H}{\partial p_r^2}=\frac{1}{m},$$ and
  \begin{multline*}
   \frac{\partial^2 H}{\partial r^2}=\frac{3}{m r^4} (p_\theta - \lambda )^2 + \frac{\mathrm{d^2} U}{\mathrm{d} r^2} = \frac{3}{r} \frac{\mathrm{d} U}{\mathrm{d} r}+ \frac{\mathrm{d^2} U}{\mathrm{d} r^2}
   \\=c\frac{3(2+\epsilon)}{r^{4+\epsilon}}-c\frac{(3+\epsilon)(2+\epsilon)}{r^{4+\epsilon}}=-c\frac{\epsilon(2+\epsilon)}{r^{4+\epsilon}}.
  \end{multline*}

  In the full system, the relative equilibrium is manifested as the unstable periodic orbits $r=r_{p_\theta}$, $p_r=0$ and $p^\pm_\theta$ such that $(p^+_\theta - \lambda )^2=(p^-_\theta - \lambda )^2$. Following general results on the persistence of normally hyperbolic invariant manifolds \cite{Fenichel71}, these periodic orbit persist if the rotational symmetry is broken, provided the perturbation is not too big. Note that according to our assumptions these periodic orbits are not associated with a local maximum of $U$.
  
  The condition $U\in o(r^{-2})$ as $r\rightarrow\infty$ is reminiscent of the assumption made by the authors of \cite{Chesnavich1980}. However, they consider a growth restriction near the origin, namely that for all $\theta$
  $$\left(\frac{\lambda^2}{2 m r^2} + U(r,\theta)\right)\in o(r^{-2}) \quad \text{as}\quad r\rightarrow0,$$
  and additionally require $\left(\frac{\lambda^2}{2 m r^2}  + U(r,\theta)\right)$ to have at most one maximum for each $\theta$. We do not impose restrictions on $U$ near $r=0$ and admit several maxima.


 \subsection{Potential energy}\label{subsec:ch4potential}

  The potential as suggested by Chesnavich \cite{Chesnavich1986} is the sum
  $$ U(r,\theta ) = U_{CH} (r) + U_{*} (r,\theta),$$
  where $U_{CH}$ is a radial long range potential and $U_{*}$ a short range ``hindered rotor'' potential that represents the anisotropy of the rigid molecule CH$_3^+$ (\cite{Jordan1991}, \cite{Chesnavich1980}). 

\begin{figure}
  \centering
    \includegraphics[width=0.45\textwidth]{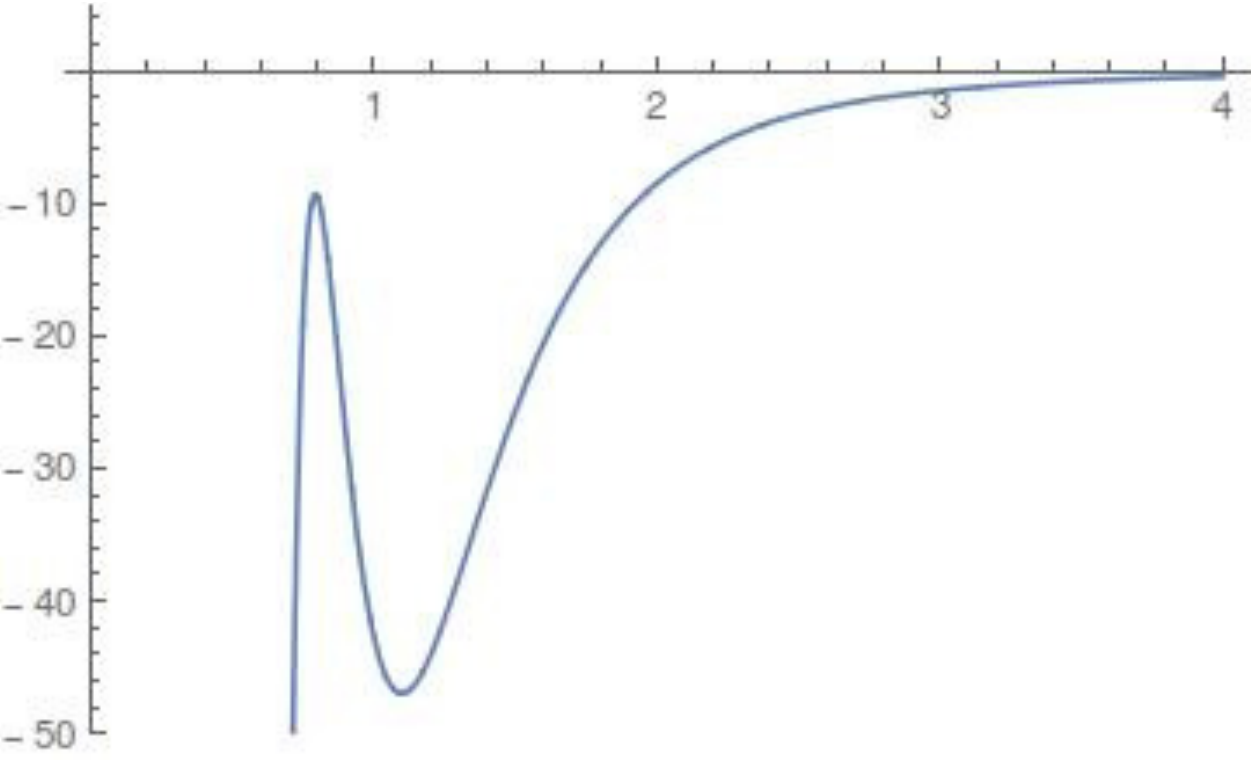}
    \includegraphics[width=0.45\textwidth]{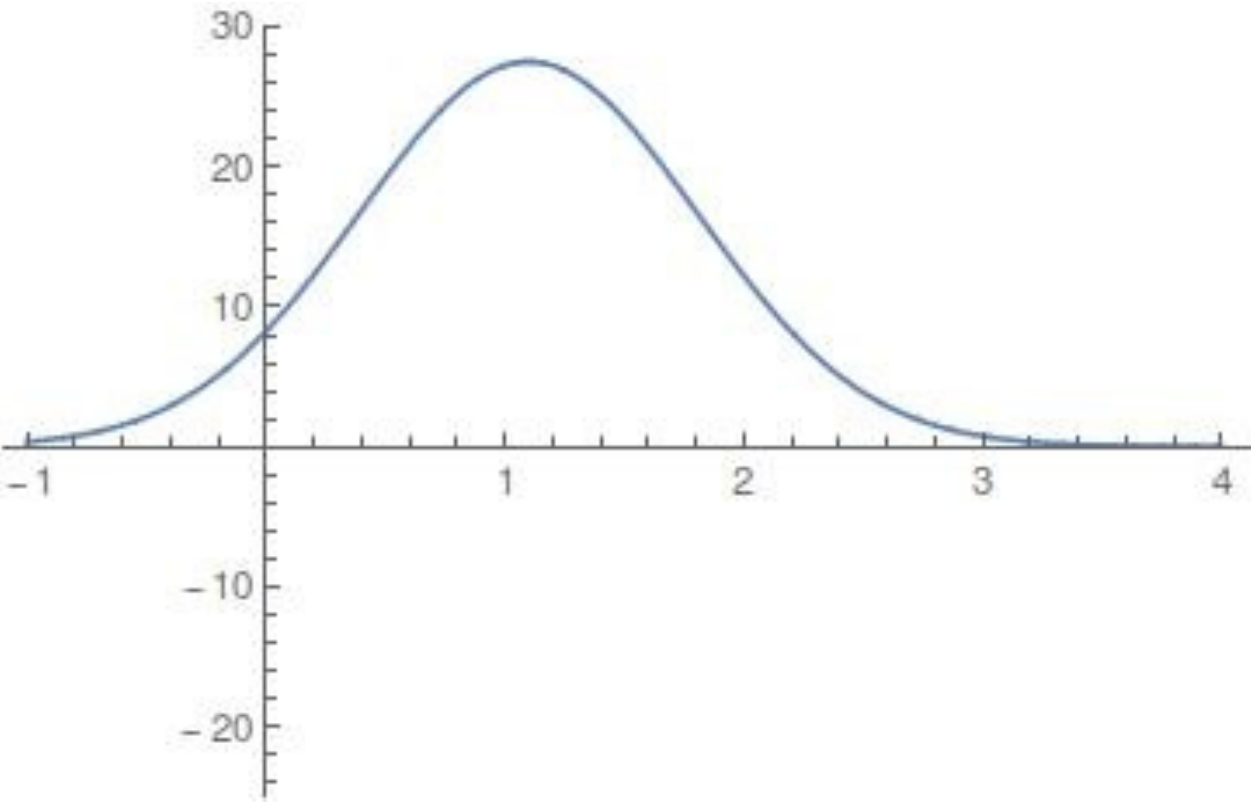}   
  \caption{Left: Graph of $U_{CH}$ versus $r$. Right: Graph of $U_{0}$ with $a = 1$ versus $r$.}
  \label{fig:ChesnavichCH4pR}
  \end{figure}

  The long range potential is defined by
  \[U_{CH} (r) =  \frac{D_e}{c_1 - 6} \left( 2 (3-c_2) e^{c_1 (1-x)}  - \left( 4 c_2 - c_1 c_2 + c_1 \right) x^{-6} - (c_1 - 6) c_2 x^{-4} \right), \]
  where $x = \frac{r}{r_e}$. The constants $D_e=47$ kcal/mol and $r_e=1.1$ \AA\ represent the C-H dissociation energy and equilibrium bond length respectively. $c_1=7.37$ and $c_2=1.61$ result in a harmonic oscillator limit with stretching frequency $3000$ cm$^{-1}$. A graph of $U_{CH}$, using Chesnavich's choice of coefficients, can be found in Figure \ref{fig:ChesnavichCH4pR}. As expected for long range interactions, it is meant to dominate the potential for large values of $r$ and not be subject to the orientation of CH$_3^+$. Therefore $U_{CH}$ is independent of the angle and its leading term for large $r$ is $r^{-4}$. Since $U_{CH}$ also dominates the short range potential in the neighbourhood of $r=0$, Chesnavich suggest a cut-off at $r=0.9$. The cut-off is not near the region of interest in our study of roaming, nor does it have any significant implications.

\begin{figure}
  \centering
  \includegraphics[width=.49\textwidth]{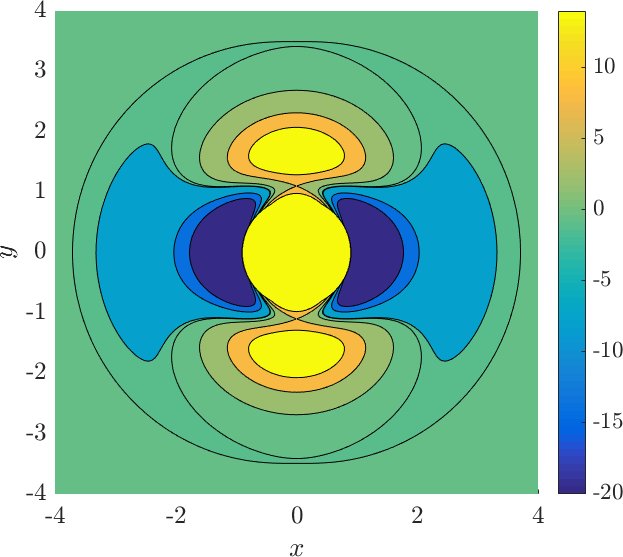}
  \includegraphics[width=.49\textwidth]{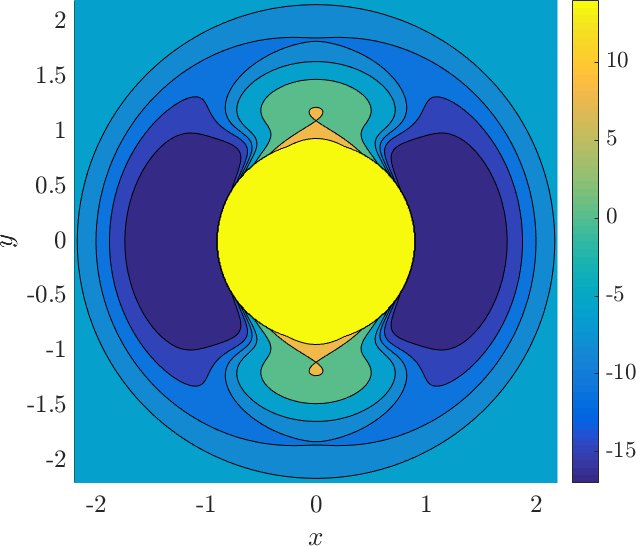}
  \caption{Contour plots of potential for $a = 1$, corresponding to a late transition, and $a = 5$, corresponding to an early transition.}
  \label{fig:potentialcontour1}
  \end{figure}

  The short range potential has the form
  $$ U_{*} (r,\theta) = \frac{U_0(r)}{2} (1 - \cos 2 \theta ),$$
  where
  $$U_0(r) = U_e e^{-a(r-r_e)^2},$$
  is the rotor barrier, which is a smoothly decreasing function of the distance $r$, and $U_e=55$ kcal/mol is the barrier height, see Figure \ref{fig:ChesnavichCH4pR}. The constant $a$ influences the value of $r$ at which the transition from vibration to rotation occurs. The transition is referred to as \emph{early} if it occurs at small $r$ and as \emph{late} otherwise. For comparison, see the late transition for $a=1$, which we will be using, and the early transition for $a=5$ in Figure \ref{fig:potentialcontour1}. Note the different proportions of the potential well (dark blue) with respect to the high potential islands along the vertical axis.
  
  Note that the angular dependence $(1 - \cos 2\theta)$ in $U_{*}$ is $\pi$-periodic and even. These properties induce a reflection symmetry of $U$ with respect to the $x$ and $y$ axes, because
  $$U(r,\theta)=U(r,-\theta),$$
  corresponds to the reflection about the $x$ axis and
  $$U(r,\theta)=U(r,-\theta+\pi),$$
  corresponds to the reflection about the $y$ axis.


 \section{Setting up the transport problem}\label{sec:dynamics}
  We follow \cite{Mauguiere2014} and set $a = 1$ for a slow transition from vibration to rotation. In what follows we also assume $\lambda=0$, unless stated otherwise. This section introduces features of the potential relevant to finding periodic orbits, defining dividing surfaces and formulating roaming in terms of transport between regions on the energy surface.

 \subsection{Energy levels and Hill regions}\label{subsec:hill}
  Here we give details about the features of the potential relevant to the dynamics of the system. Being the most basic characteristic of the potential, we look at critical points of the potential that give valuable information about local dynamics and at level sets that tell us about the accessible area in configuration space.

  Due to the reflection symmetry of the potential about the $x$ and $y$ axes introduced above, critical points always come in pairs. We will denote them by $q_i^\pm$, where $i$ indicates the index of the critical point and the superscript $+$ stands for the upper half plane $\theta\in[0,\pi)$, while $-$ stands for the lower. Here we present a list of critical points:
  \begin{itemize}
  \item $q_0^\pm$ - two wells at $(r,\theta)=(1.1, 0)$ and $(1.1, \pi)$ with $U(q_0^\pm)=E_0\approx -47$,
  \item $q_1^\pm$ - two index-$1$ saddles at $(3.45, \frac{\pi}{2})$ and $(3.45, \frac{3\pi}{2})$ with $U(q_1^\pm)=E_1\approx -0.63$,
  \item $\widetilde{q}_1^\pm$ - two index-$1$ saddles at $(1.1, \frac{\pi}{2})$ and $(1.1, \frac{3\pi}{2})$ with $U(\widetilde q_1^\pm)=\widetilde{E}_1\approx 8$,
  \item $q_2^\pm$ - two index-$2$ saddles at $(1.63, \frac{\pi}{2})$ and $(1.63, \frac{3\pi}{2})$ with $U(q_2^\pm)=E_2\approx 22.27$.
  \end{itemize}

  The potential wells correspond to the two isomers of CH$_4^+$ with the free H atom close to the CH$_3^+$ molecule. All index-$1$ saddles are involved in isomerisation and the two index-$2$ saddles provide us with interesting geometries of the accessible regions in configuration space. For zero angular momentum ($\lambda = 0$), the critical phase space points of $H$ are given by $ z_i^\pm = ( q_i^\pm, 0)$ and $ \widetilde{z}_1^\pm = (\widetilde{q}_1^\pm, 0)$.

  The critical energies are ordered as 
  $$E_0 < E_1 < 0 < \widetilde{E}_1 < E_2,$$
  and all critical points can be found in the contour plot in Figure \ref{fig:potentialcontour1}.
  
  For a given fixed energy $E$, we are interested in the accessible region in the configuration space which, following the celestial mechanics literature, we refer to as Hill region \cite{Hill1905}, and the geometry of the energy surface. Since the system as defined in Section \ref{subsec:system} is a natural mechanical system $H=T+U$ and the kinetic energy $T$ is always non-negative, the Hill region consists of all $(r,\theta)$ such that $U(r,\theta) \leq E$ and is bounded by the equipotential $U(r,\theta) = E$.
  
  To see what the Hill regions look like, we note that the two wells $q_0^\pm$ give rise to two topological discs that connect into an annulus at $E=E_1$ via $q_1^\pm$. With $E\rightarrow0$ the annulus widens until at $E=0$ it loses compactness and covers the whole plane except for a disc near the origin. This cut-out disc decomposes at $E=\widetilde{E}_1$ into three, two areas of high potential around $q_2^\pm$ and the cut-off of the potential at $r=0.9$ mentioned earlier. Above $E=E_2$ only the cut-off at $r=0.9$ remains inaccessible. Topologically this is equivalent to the case with the energy $0<E<\widetilde{E}_1$.
  
  Hill regions are shown for various energies in Figure \ref{fig:potentialcontour1}. For comparison, we also include Hill regions for the case $a=5$ in Figure \ref{fig:potentialcontour1}, where the transition from vibration to rotation occurs earlier. Although energy levels remain topologically equivalent, note the larger potential well and the smaller energy interval where the boundary of Hill region consists of three circles.


 \subsection{Relevant periodic orbits}\label{subsec:PO}

  Next we study the invariant structures that can be found in the system at various energies. Critical points $z_i^\pm$ described in Section \ref{subsec:hill} are the most basic invariant structures at energies $E_i$. In the following we discuss (non-degenerate) periodic orbits on the $3$-dimensional energy surface. They create a backbone for the understanding the dynamical behaviour of our system.
 
  Similarly to critical points, periodic orbits also come in pairs because of the symmetry of the potential. The periodic orbits are then related by the discrete rotational symmetry $$(r,\theta,p_r,p_\theta) \mapsto(r,\theta+\pi,p_r,p_\theta),$$ or the discrete reflection symmetry $$(r,\theta,p_r,p_\theta) \mapsto(r,-\theta,p_r,-p_\theta).$$ In contrast to critical points, non-degenerate periodic orbits persist in energy intervals forming one-parameter families. As periodic orbits evolve with varying energy, they occasionally bifurcate with other families of periodic orbits.
  
  Based on the knowledge of Hill regions we gained in Section \ref{subsec:hill}, we can formulate some expectations about periodic orbits in this system.
  For $E\leq E_1$, the system does not admit rotating periodic orbits, orbits that are periodic in $\theta$ and along which always $p_\theta>0$ or $p_\theta<0$. Rotating orbits project onto the configuration space as circles with the origin contained in their interior. Instead in the interval $E_0<E\leq E_1$ we can only expect vibrating, oscillator like, periodic orbits. Of special interest are periodic orbits that project onto a line with both ends on an equipotential. In the celestial mechanics literature these orbits are referred to as periodic \emph{brake orbits} (the name is due to Ruiz \cite{Ruiz75}).
  
\begin{figure}
  \centering
  \includegraphics{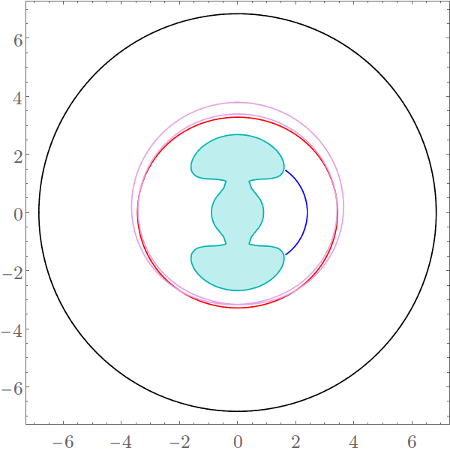}
  \caption{Configuration space projections of $\Gamma^i_\pm$ (blue), $\Gamma^o_\pm$ (black), $\Gamma^a_\pm$ (red) and one orbit of the family $\Gamma^b$ (magenta) at energy $E=2$.}
  \label{fig:po2}
  \end{figure}

  Here we present a list of the important families of periodic orbits together with a brief description of their evolution. Configuration space projections of the periodic orbits at $E=2$ are shown in 
  Fig.~\ref{fig:po2}.

\begin{figure}
  \centering
  \includegraphics[width=10cm]{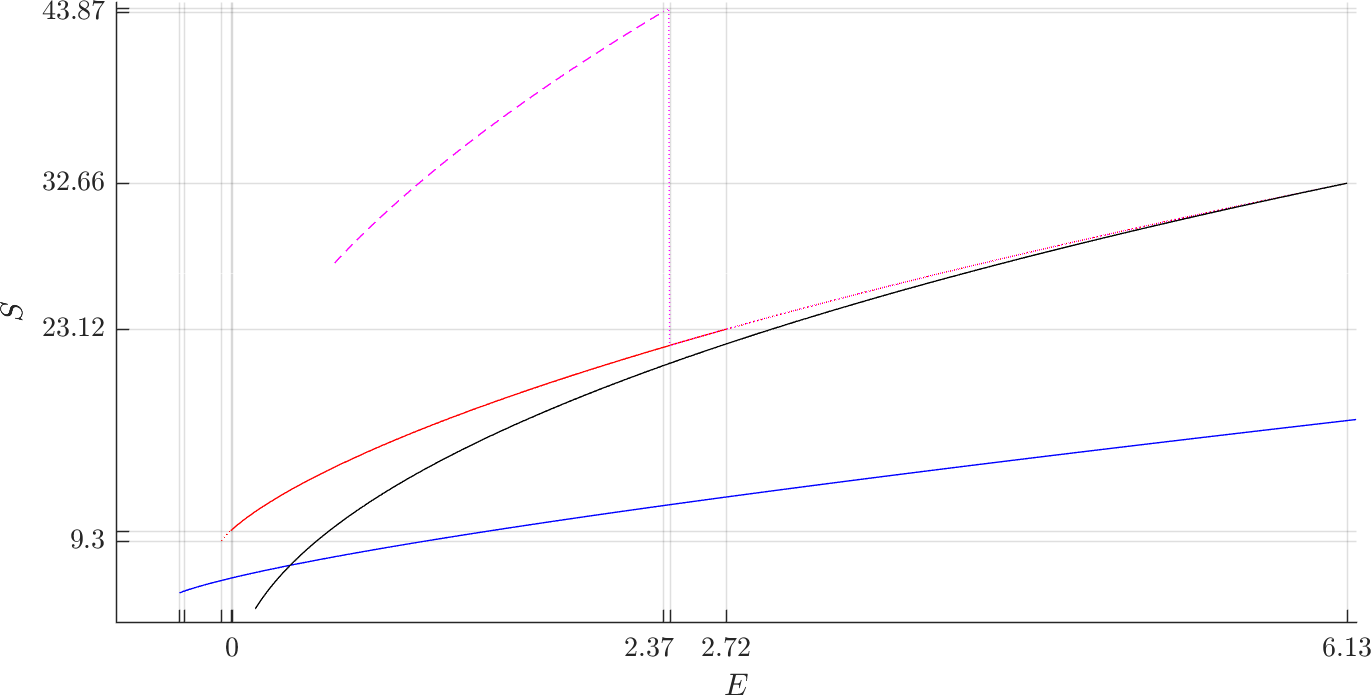}
  \caption{Bifurcation diagrams showing $\Gamma^i_\pm$ (blue), $\Gamma^o_\pm$ (black), $\Gamma^a_\pm$ (red) and orbits of the family $\Gamma^b$ (magenta) in the energy-action ($E,S$) plane.}
  \label{fig:bifurcationS}
  \end{figure}

  \begin{itemize}
  \item $\Gamma^i$: The family of periodic orbits $\Gamma^i$ is born in a saddle-centre bifurcation at energy $E = -.29$. Until a host of bifurcations above $E=20$,  $\Gamma^i$ consists of hyperbolic brake orbits. Around $E=21.47$ the orbits become inverse hyperbolic and at $E = 22.27$ they become heteroclinic to $z_2^\pm$ and undergo a Morse bifurcation, similar to those described in \cite{MacKay2014}. At higher energies, $\Gamma^i$ consists of rotating orbits that undergo further bifurcations. The periodic orbits are by some authors referred to as inner or tight periodic orbits. We denote the individual orbits by $\Gamma^i_+$ and $\Gamma^i_-$. For $E \leq 22.27$ $\Gamma^i_\pm$ is the brake orbit in the potential well associated with $z^\pm_0$ and for $E > 22.27$ the subscript $\pm$ corresponds to the sign of $p_\theta$ along the rotating periodic orbit.

  \item $\Gamma^o$: This family of unstable periodic orbits originates at $r=\infty$ at $E=0$. With increasing energy the orbits monotonously decrease in radius and remain unstable until a bifurcation with $\Gamma^a$ and $\Gamma^b$ at $E=6.13$, where $\Gamma^a$ and $\Gamma^b$ are described below. These periodic orbits are sometimes called outer or orbiting periodic orbits, because these are the periodic orbits with the largest radius at the energies where they exist. We denote the individual orbits with $p_\theta>0$ and $p_\theta<0$ by $\Gamma^o_+$ and $\Gamma^o_-$ respectively.

  \item $\Gamma^a$: These periodic orbits are created in a saddle-centre bifurcation at $E=-.0602$ as stable, turn unstable at $E=-.009$ and remain unstable until a period doubling bifurcation at $2.72$. The family disappears in the aforementioned bifurcation with $\Gamma^o$ and $\Gamma^b$. At all energies, the configuration space projection of $\Gamma^a$ is located between those of $\Gamma^i$ and $\Gamma^o$ and we will refer to the orbits as the middle periodic orbits. We denote the individual orbits with $p_\theta>0$ and $p_\theta<0$ by $\Gamma^a_+$ and $\Gamma^a_-$ respectively.

  \item $\Gamma^b$: The product of a saddle-centre bifurcation at $E=-.0023$ that quickly becomes inverse hyperbolic. Around $E=2.37$ the orbits become elliptic and undergo a reverse period doubling at $E=2.4025$. Note that the energetic gap between these two bifurcations is so small that they are almost indistinguishable in Figure \ref{fig:bifurcationS}. After that $\Gamma^b$ remains stable until it collides with $\Gamma^a$ and $\Gamma^o$. For $E<2.4025$ the family consists of four periodic orbits with twice the period compared to all the previously mentioned ones. The orbits related by discrete symmetries mentioned above.
  \end{itemize}

  $\Gamma^i$ is important because its orbits lie in the potential well and have the largest radial coordinate $r$ of all periodic orbits in the well. $\Gamma^o$ are the outermost periodic orbits and trajectories with a larger radial coordinate $r$ and $p_r>0$ go to infinity in forward time, i.e. $r\rightarrow\infty$ as $t\rightarrow\infty$. As mentioned above, the configuration space projections of $\Gamma^a$ lie between $\Gamma^i$ and $\Gamma^o$. In fact there are no other periodic orbits with single period ($2\pi$-periodic in $\theta$) in this region of configuration space. We use orbits of the family $\Gamma^a$ in Section \ref{subsec:divide energy surface} to define dividing surfaces and divide phase space into regions. $\Gamma^b$ is needed for a complete description of the evolution of $\Gamma^o$ and $\Gamma^a$ and its bifurcations may hint at qualitative changes of structures formed by invariant manifolds. 
  
  There are various other periodic orbits, most notably ones corresponding to stable vibrations of the bound CH$_4^+$ molecule, Lyapunov orbits associated with $z_1^\pm$ and $\widetilde {z}_1^\pm$ that play a role in isomerisation and periodic orbits involved in various bifurcations with the orbits mentioned above. All of these will not play a role in our further considerations.

  With non-zero angular momentum, periodic orbits of a family remain related by the discrete rotational symmetry, but not by the discrete reflection symmetry and some other properties are different too. The inner periodic orbits are no longer brake orbits for $\lambda\neq 0$ and their projections onto configuration space are topological circles instead of lines. Similarly rotating orbits of the same family do not have the same configuration space projection and bifurcate at different energies. With increasing $|\lambda|$ the differences become more pronounced.

  In Figure \ref{fig:bifurcationS} we present the evolution of the orbits in above-mentioned families in the energy-action ($E,S$) plane. We will explain in Section \ref{subsec:DS} why flux through a dividing surface associated with a vibrating periodic orbit is equal to its action, while for rotating periodic orbits it is equal to twice its action.
  
\begin{figure}
  \centering
  \includegraphics[width=10cm]{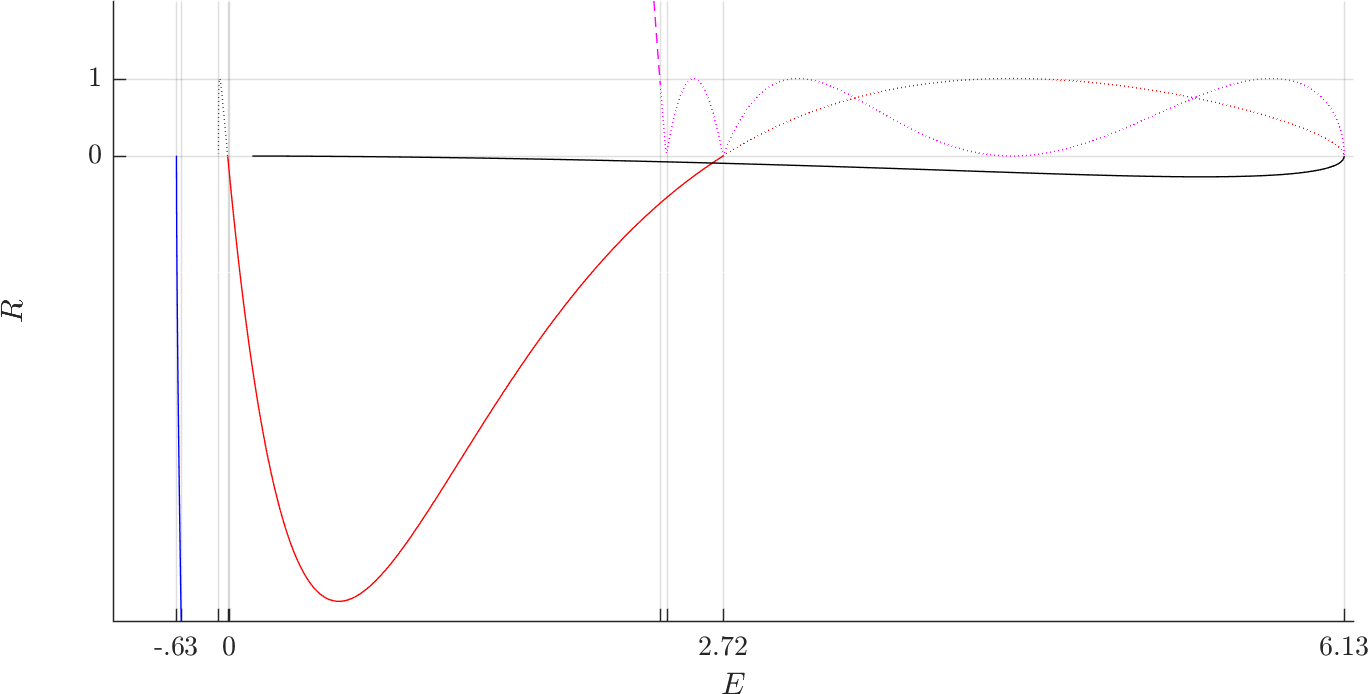}
  \caption{Bifurcation diagrams showing $\Gamma^i_\pm$ (blue), $\Gamma^o_\pm$ (black), $\Gamma^a_\pm$ (red) and orbits of the family $\Gamma^b$ (magenta) in the energy-residue ($E,R$) plane.}
  \label{fig:bifurcationR}
  \end{figure}

  Figure \ref{fig:bifurcationR} shows the evolution of the Greene residue of orbits in the families. The Greene residue, due to J. M. Greene \cite{Greene68} is a quantity characterizing the stability of the orbits. It is derived from the monodromy matrix, a matrix that describes the behaviour of solutions in the neighbourhood of a periodic orbit.

  The monodromy matrix and the Greene residue are defined as follows. For a periodic orbit $\Gamma$ with the parametrisation $\gamma(t)$ and period $T$, let $M(t)$ be the matrix satisfying the variational equation
  \begin{equation*}
  \dot{M}(t)=JD^2 H(\gamma(t))M(t),
  \end{equation*}
  where $J=\begin{pmatrix} 0 & Id\\ -Id & 0\end{pmatrix}$,
  with the initial condition
  $$M(0)=Id.$$
  The monodromy matrix is defined by $M=M(T)$. It describes how an initial deviation $\delta$ from $\gamma(0)$ changes after a full period $T$. For $\delta$ sufficiently small the relationship is
  \begin{equation*}
  \Phi_H^T(\gamma(0)+\delta)=\gamma(T)+M\delta+O(\delta^2),
  \end{equation*}
  where $\Phi_H^t$ is the Hamiltonian flow.

  If $\delta$ is an initial displacement along the periodic orbit $\delta\parallel J\nabla H$, then by periodicity $\delta$ is preserved after a full period $T$, i.e. $M\delta=\delta$. A similar argument holds for an initial displacement perpendicular to the energy surface $\delta\parallel \nabla H$. Consequently, two of the eigenvalues of $M$ are $\lambda_1=\lambda_2=1$. More details including a reduction of $M$ can be found in \cite{Eckhardt91}.

  As the variational equation satisfied by $M(t)$ is Hamiltonian, the preservation of phase space volume following Liouville's theorem implies that the determinant $\det M(t)=\det M(0)=1$ for all $t$. Therefore for the two remaining eigenvalues we have $\lambda_3\lambda_4=1$ and we can write them as $\lambda$ and $\frac{1}{\lambda}$. $\Gamma$ is hyperbolic if $\lambda>1$, it is elliptic if $|\lambda|=1$ and it is inverse hyperbolic if $\lambda<-1$.

  \begin{definition}
  The Greene residue of $\Gamma$ is defined as
  $$R=\frac{1}{4}(4-Tr M),$$
  where $M$ is the monodromy matrix corresponding to the periodic orbit $\Gamma$. 
  \end{definition}

  Knowing that $\lambda_1=\lambda_2=1$, we can write $R$ as
  $$R=\frac{1}{4}\left(2-\lambda-\frac{1}{\lambda}\right).$$
  By definition $R<0$ if $\Gamma$ is hyperbolic, $0<R<1$ if it is elliptic and $R>1$ if it is inverse hyperbolic.


 \subsection{Transition states and dividing surfaces}\label{subsec:DS}

  In this section we discuss dividing surfaces associated with transition states, the backbone of Transition State Theory. Following \cite{MacKay2014,MacKay2015} we define transition states more formally as follows.
  
  \begin{definition}[TS]
    A transition state for a Hamiltonian system is a closed, invariant, oriented, codimension-$2$ submanifold of the energy surface that can be spanned by two surfaces of unidirectional flux, whose union divides the energy surface into two components and has no local recrossings.
  \end{definition}
  
  The name transition state is due to the fact it is a structure found between areas of qualitatively different types of motion, a transition between two types of motion so to say. One can imagine the transition between types of motion corresponding to physical states like reactants and products or the transition between rotation and vibration.
  
  In a system with $2$ degrees of freedom, a TS consists of unstable periodic orbit. Generally a TS is a codimension-$2$ normally hyperbolic invariant manifold, a manifold on the energy surface invariant under the Hamiltonian flow, such that instabilities transversal to it dominate the instabilities tangential to it (\cite{Fenichel71}, \cite{Hirsch77}).

  In general, a dividing surface (DS) is a surface that divides the energy surface into two disjoint components. By a DS associated with a TS we mean a union of the two surfaces of unidirectional flux that is constructed as follows.
  
  For a fixed energy $E$, let $(r_\Gamma,\theta_\Gamma)$ be the projection of the periodic orbit $\Gamma$ onto configuration space, then the DS is the surface $(r_\Gamma,\theta_\Gamma,p_r,p_\theta)$, where $(p_r,p_\theta)$ are given implicitly by the energy equation
  $$E= \frac{1}{2 m} p_r^2 + \frac{1}{2 I} p_\theta^2  + \frac{1}{2 m r^2} (p_\theta - \lambda)^2+U(r,\theta).$$
  This construction also works for stable periodic orbits, but the resulting DS admits local recrossings.
  In the following a DS associated with a TS is always the surface constructed this way.

  We will refer to the DSs associated with $\Gamma^i$, $\Gamma^o$ and $\Gamma^a$ as inner, outer and middle, respectively. For our investigation, we do not need to distinguish between the DSs associated to $\Gamma^i_+$ and $\Gamma^i_-$, therefore we always refer to the former unless explicitly stated otherwise. We are mainly interested in the influence of local energy surface geometry on the geometry of DSs and in the dynamics on DSs under the corresponding return map.

  The geometry of the DSs is due to the form of the kinetic energy and the local geometry of the energy surface. It is well known that a DS associated to a brake periodic orbit is a sphere and the brake periodic orbit is an equator of this sphere, \cite{Waalkens04}. The equator divides the sphere into hemispheres, whereby the flux through the two hemispheres is equal in size and opposite in direction. Trajectories passing this sphere from reactants to products intersect one hemisphere and the other hemisphere is crossed on the way from products to reactants. The flux through a hemisphere is then by Stokes' theorem equal to the action of the periodic orbit \cite{MacKay1990}.

  Rotating periodic orbits, on the other hand, such as $\Gamma^o$, give rise to a DS that is a torus. The two orbits of the same family with opposite orientation are circles on the torus and divide it into two annuli with properties identical to the hemispheres. Using Stokes' theorem we find that the flux across each annulus is given by the sum of the actions of the two orbits, or simply twice the action of a single orbit \cite{MacKay2014}.
  
  Should it be necessary to distinguish the hemispheres or annuli of a DS by the direction of flux, the outward hemisphere or annulus is the one intersected by the prototypical dissociating trajectory defined by $\theta=0$, $p_r>0$, $p_\theta=0$ and/or $\theta=\pi$, $p_r>0$, $p_\theta=0$. The inward hemisphere or annulus is then intersected by $\theta=0$, $p_r<0$, $p_\theta=0$ and/or $\theta=\pi$, $p_r<0$, $p_\theta=0$


 \subsection{Division of energy surface}\label{subsec:divide energy surface}

  Using the inner and outer DSs we can define regions on the energy surface and formulate roaming as a transport problem.

  The area bounded by the surface $r=0.9$ and the two inner DSs represents the two isomers of CH$_4^+$. We denote the two regions by $B_1^+$ and $B_1^-$. The unbounded region beyond the outer DS, denoted $B_3$, represents the dissociated molecule.

  It is therefore in the interaction region between the inner and the outer DS, denoted $B_2$, where the transition between CH$_4^+$ and CH$_3^++$H occurs. When in $B_2$, the H atom is no longer in the proximity of CH$_3^+$, but still bound to the CH$_3^+$ core. This is the region, where the system exhibits roaming. Contained in $B_2$ are $\Gamma^a$ and various other periodic orbits that may play a role in roaming.

  Dissociation can in this context be formulated as a problem of transport of energy surface volume from $B_1$ to $B_3$. Such volume contains trajectories that originate in the potential well, pass through the interaction region and never return after crossing the outer DS. Since each trajectory passing from $B_1$ to $B_2$ crosses the inner DS and leaves $B_2$ by crossing the outer DS, we may restrict the problem to the interaction region. Because roaming is a particular form of dissociation, it too has to be subject to transport from the inner DS to the outer DS.

  It is well known that transport to and from a neighbourhood of a unstable periodic orbits is governed by its stable and unstable invariant manifolds. The problem can be reformulated accordingly. This means, of course, by studying the structure of heteroclinic intersections of stable and unstable invariant manifolds of $\Gamma^i$ and $\Gamma^o$, as well as with $\Gamma^a$ that, as we will soon see, sits inside the homoclinic tangle of $\Gamma^o$.

  We will denote the invariant manifolds of $\Gamma^i_+$ by $W_{\Gamma^i_+}$. We will further use a superscript $s$ and $u$ to label the stable and unstable invariant manifolds and add an extra superscript $-$ and $+$ for the branches that leave the neighbourhood of $\Gamma^i_+$ to the CH$_4^+$ side ($r$ smaller) or to the CH$_3^++$H side ($r$ larger), respectively. $W_{\Gamma^i_+}^{u+}$ therefore denotes the unstable branch of the invariant manifolds of $\Gamma^i_+$ that leaves the neighbourhood of $\Gamma^i_+$ to the CH$_3^++$H side. Invariant manifolds of other TSs will be denoted analogously.

  We remark that we may use TST to consider the evolution of periodic orbits in the energy-action plane shown in Figure \ref{fig:bifurcationS} in the context of transport of energy surface volume from $B_1$ to $B_3$. Recall for Section \ref{subsec:DS} that the flux across the outer and middle DSs is twice the action of $\Gamma^o_+$ and $\Gamma^a_+$ respectively. The combined flux through both inner DSs is twice the action of $\Gamma^i_+$.  We see that for $E\leq.32$, the outer DS has the lowest flux, while for higher energies it is the inner DS.


 \section{Dynamics of the Chesnavich model}\label{sec:observations}

  Before we proceed to the discussion of how invariant manifolds cause slow dissociations, let us describe some numerical observations of how the system behaves in certain phase space regions. The observations will later be explained using invariant manifolds. In the following, we offer insight into the amount of time needed to dissociate, the locations where dissociation is fast or slow and how these properties change with increasing energy. We use this knowledge to establish a link between invariant manifolds and slow dissociation on which we further elaborate in Section \ref{sec:discussion} in the context of roaming.

 \subsection{Residence times and rotation numbers}\label{subsec:res th0}

  For various energies $0<E<6.13$ where $\Gamma^o$ exists, we investigate trajectories starting in $B_1^+$, $B_2$ and $B_3$ on the surface $\theta=0, p_\theta>0$. We study how long it takes trajectories to reach a terminal condition representing the dissociated state.

  In Section \ref{subsec:divide energy surface} we said that we consider the molecule dissociated as soon as the system enters $B_3$. Naturally, then the terminal condition should be that trajectories reach the outer DS. However using the outer DS raises uncertainty of whether a faster dissociation is a dynamical property or a result of the changing position of $\Gamma^o_\pm$ with energy. To prevent this uncertainty, we use a fixed terminal condition. Since for $E\rightarrow0$, the radius of $\Gamma^o_\pm$ diverges, no fixed terminal condition can represent the dissociated state for all energies. We decided to define the terminal condition by $r=15$ that works well for $E\geq0.4$ at the cost of losing the energy interval $E<0.4$.

  In the following we consider residence times and rotation numbers, i.e. time and change in angle needed for trajectories starting on $\theta=0, p_\theta>0$ to reach the surface $r=15$ in $B_3$.

  \begin{figure}
  \includegraphics[width=.5\textwidth]{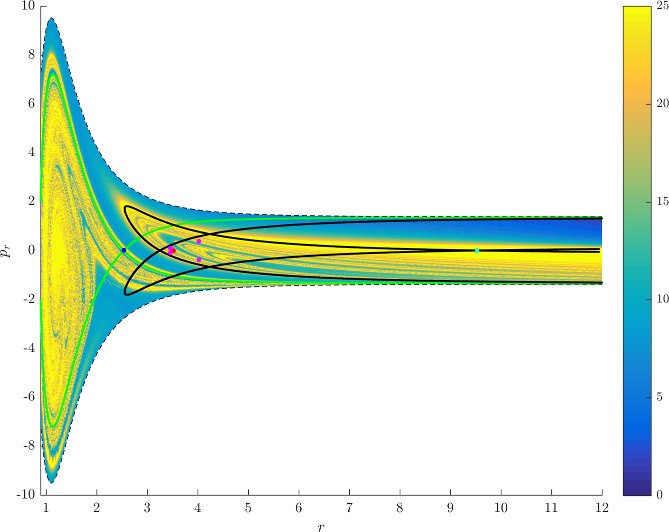}
  \includegraphics[width=.5\textwidth]{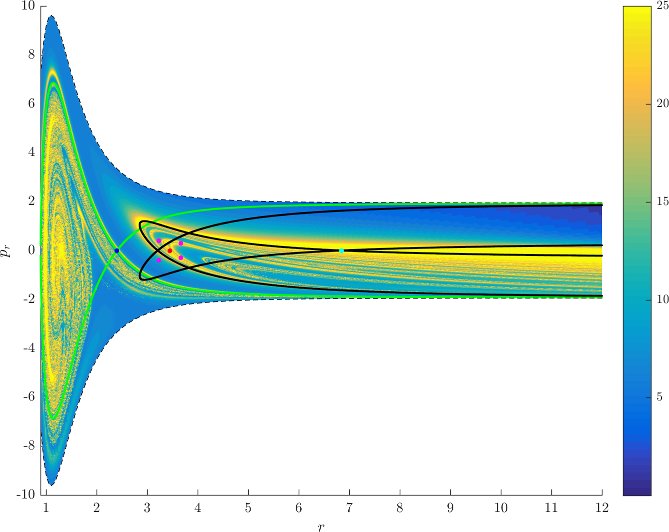}\\
  \includegraphics[width=.5\textwidth]{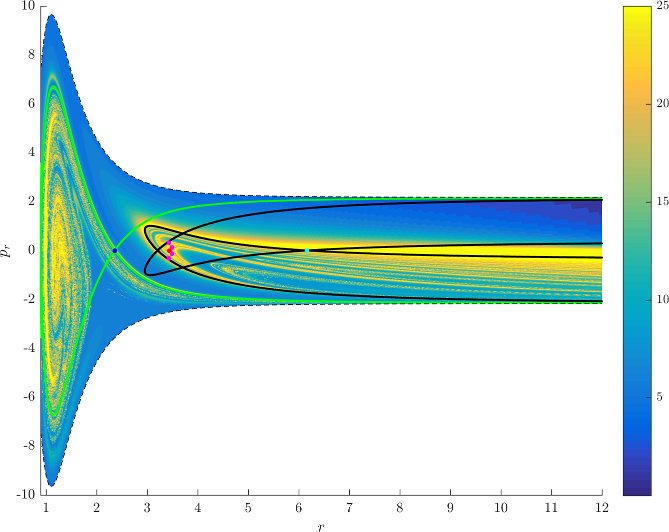}
  \includegraphics[width=.5\textwidth]{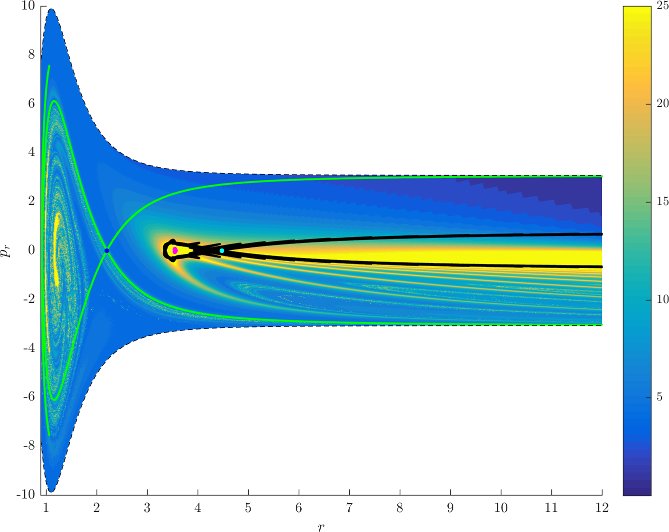}
  \caption{From top left to bottom right the plots show residence times on the surface of section $\theta=0$, $p_\theta>0$ for energies $E=1,2,2.5,5$. The dots correspond to the periodic orbits $\Gamma^i_+$ (blue), $\Gamma^a_+$ (red), orbits of the family $\Gamma^b$ (magenta) and $\Gamma^o_+$ (cyan). Invariant manifolds of $\Gamma^i_+$ (green) and $\Gamma^o_+$ (black) are also included.}
  \label{fig:rot}
  \end{figure}
  
  Figure \ref{fig:rot} shows rotation numbers for selected energies, with marked periodic orbits and invariant manifolds. As expected, initial conditions with $p_r>0$ large are the fastest ones to escape. The slowest ones are located near the periodic orbits and near $p_r=0$ ($p_\theta$ large).

  For $E\leq 2.5$, almost all initial conditions in $B_1^+$ were slow to escape. For higher energies, most of the slow dissociation occurs around $\Gamma^a_+$ and $\Gamma^o_+$, the slowly dissociating trajectories have a negative initial $p_r$ very close to zero. This observation is easily explained by noting that configuration space projections of these trajectories are almost circular and spend most of the time in the region where the potential is very flat and almost independent of $\theta$, thus $\dot {p}_\theta\approx0$.

  Chaotic structures that can be seen in $B_1^+$ are the result of lengthy escape from a potential well. The only known structure responsible for fractal-like patterns and one closely linked to chaotic dynamics are invariant manifolds, in this case $W_{\Gamma^i_+}$. Note that at $E=5$, it seems that $W_{\Gamma^o_+}$ slows the dynamics down considerably more than $W_{\Gamma^i_+}$.
  
  \begin{figure}
  \includegraphics[width=.5\textwidth]{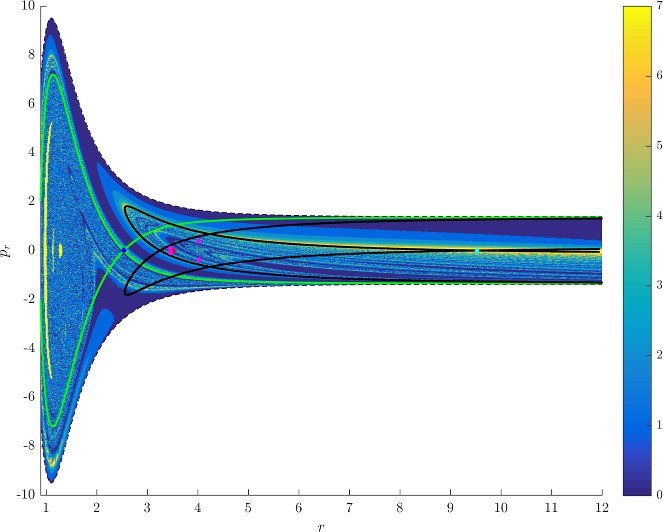}
  \includegraphics[width=.5\textwidth]{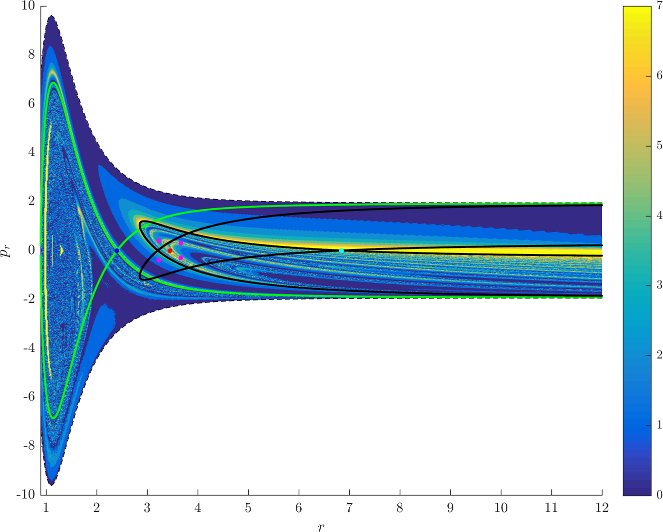}\\
  \includegraphics[width=.5\textwidth]{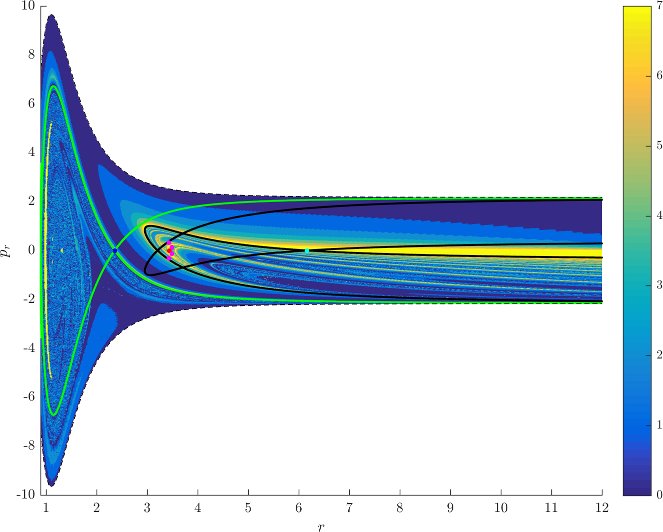}
  \includegraphics[width=.5\textwidth]{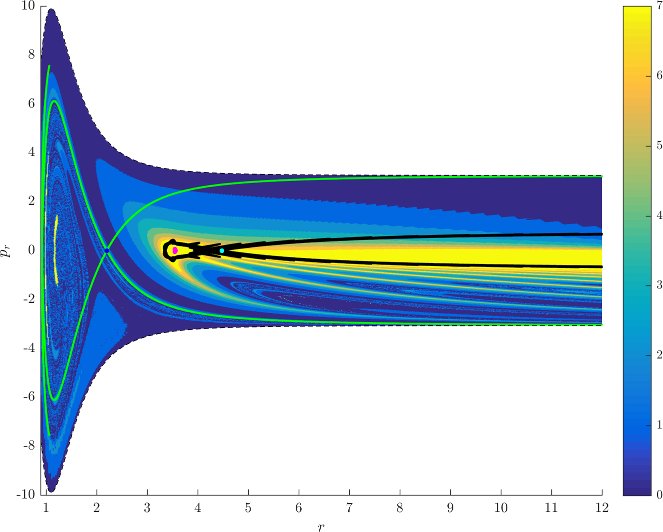}
  \caption{From top left to bottom right the plots show rotation numbers on the surface of section $\theta=0$, $p_\theta>0$ for energies $E=1,2,2.5,5$. The dots correspond to the periodic orbits $\Gamma^i_+$ (blue), $\Gamma^a_+$ (red), orbits of the family $\Gamma^b$ (magenta) and $\Gamma^o_+$ (cyan). Invariant manifolds of $\Gamma^i_+$ (green) and $\Gamma^o_+$ (black) are also included.}
  \label{fig:rot_num}
  \end{figure}

  Rotation numbers, i.e. number of completed full rotations upon dissociation, closely match residence times suggesting that slowly dissociating trajectories are ones that rotate in $B_2$ and $B_3$ for a long time. More pronounced, due to the discrete nature of the number of rotations, are structures inside $B_1^+$, just below $p_r=0$ and in the neighbourhood of $\Gamma^o_+$ and $W_{\Gamma^o_+}$. 

  Note in Figure \ref{fig:rot_num} that the fractal like structures recede with increasing energy and by $E=5$ most of them lie either in $B_1^+$, near $p_r=0$ as mentioned above and in the proximity of the homoclinic tangle of $\Gamma^o_+$. The homoclinic tangle seems to tend to a homoclinic loop as it disappears for $E \rightarrow 6.13$. It is also worth noting that fast and simple dissociation, i.e. low residence time and low rotation number, is not only becoming more dominant, but also speeding up, see Figures \ref{fig:rot} and \ref{fig:rot_num}. Due to the increase in kinetic energy in the angular degree of freedom, the dissociating trajectories are naturally not becoming more direct with increasing energy.

 \subsection{Residence times on the inner DS}\label{subsec:res inner DS}

  Similarly to the surface $\theta=0, p_\theta>0$, we can study residence times and rotation numbers for trajectories starting on a DS. In Section \ref{subsec:divide energy surface} we formulated our problem as a transport problem from the inner to the outer DS. According to Section \ref{subsec:DS}, trajectories enter $B_1^+$ through one hemisphere of the inner DS and leave through the other. Naturally we are interested in the latter hemisphere.

  \begin{figure}
  \centering
  \includegraphics[width=.49\textwidth]{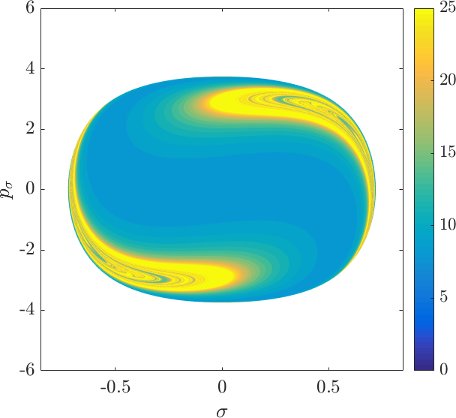}
  \includegraphics[width=.49\textwidth]{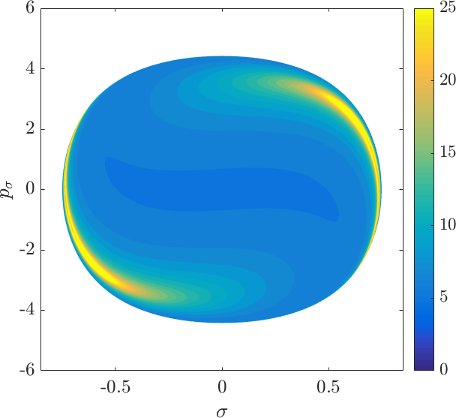}\\
  \includegraphics[width=.49\textwidth]{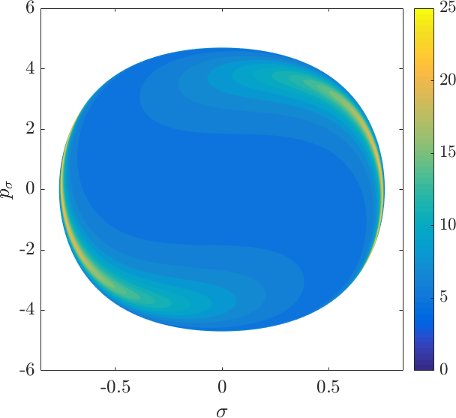}
  \includegraphics[width=.49\textwidth]{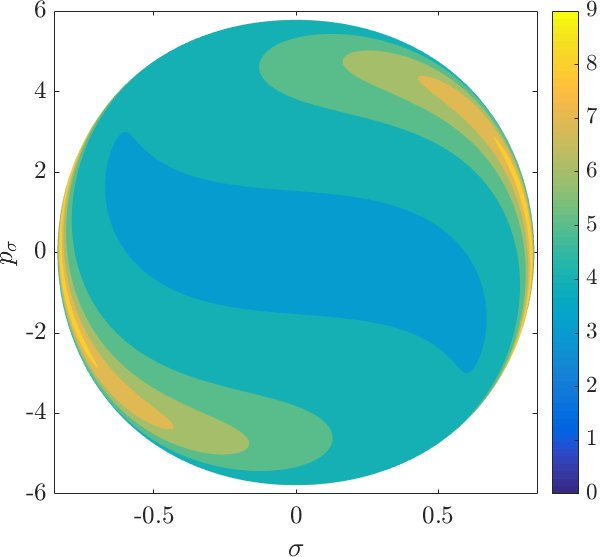}
  \caption{Residence times for initial conditions on the inner DS with outward direction for energies $E=1,2,2.5,5$. Note that the scale for $E=5$ is different, because $9$ is an upper bound for the residence time for initial conditions on the inner DS.}
  \label{fig:disI_rot}
  \end{figure}
  
  Although it is not absolutely indispensable for qualitative purposes, we prefer to work on the DS in canonical coordinates. Due to the preservation of the canonical $2$-form by the Hamiltonian flow, if we use canonical coordinates, the map from one surface of section to another is area preserving. Consequently areas of initial conditions on the inner DS corresponding slow or fast dissociation can be directly compared to the areas on the surface of section $\theta=0$, $p_\theta>0$.

  Canonical coordinates are obtained by defining a new radial variable $\rho(r,\theta)=r-\bar{r}(\theta)$ that is constant along $\Gamma^i_+$, where the curve $\bar{r}(\theta)$ is the approximation of the configuration space projection of $\Gamma^i_+$, similarly to \cite{Jaffe1999}. Due to the symmetry of the system, $\Gamma^i_+$ can be very well approximated by a quadratic polynomial for every energy. Next we use the generating function (type 2 in \cite{Arnold76})
  $$G(r,\theta,p_\rho,p_\sigma)=(r-\bar{r}(\theta))p_\rho+\theta p_\sigma.$$
  From that we obtain
  \begin{eqnarray*}
  p_r=\frac{\partial G}{\partial r}&=&p_\rho,\\
  p_\theta=\frac{\partial G}{\partial \theta}&=&p_\sigma-\bar{r}'(\theta) p_\rho,\\
  \end{eqnarray*}
  and therefore $p_\sigma=p_\theta+\bar{r}'(\theta) p_\rho$.
  The surface of section is now defined by $\rho=0$, $\dot{\rho}>0$, i.e. the outward hemisphere of the inner DS corresponding to transport in the direction from $B_1^+$ to $B_2$.
  
  Figure \ref{fig:disI_rot} shows the distribution of residence times for initial conditions on the inner DS. We can see that slow dissociation is specific to two areas of the surface of section. Initial conditions on the rest of the surface leave $B_2$ quickly. Information from the two surfaces of section suggests that $W_{\Gamma^o_+}$, and eventually $W_{\Gamma^a_+}$, intersect the inner DS in the area with long dissociation times. The areas of slow dissociation are the most pronounced for low energies, at $E=2.5$ they almost disappear. At $E=5$ we see no sign of slow dissociation, the longest residence time found at the current resolution ($6000\times6000$ initial conditions) was below $9$. This suggests that the structure responsible for roaming disappears at an energy below $2.5$. Note that even the slowest dissociation at $E=5$ takes as long as the fastest ones at $E=1$ or $E=2$.

  In summary we can say that the system exhibits various types of dissociation ranging from fast and direct, where the H atom escapes almost radially, to slow that involves H revolving a multitude of times around CH$_3^+$. Long dissociations seem to occur in fractal-like structures that are caused by invariant manifolds, proof of which will be given in Section \ref{subsec:sec manifs}.


 \subsection{Sections of manifolds}\label{subsec:sec manifs}

  Let us now have a closer look at manifolds on the two surfaces of section presented above and establish a link between invariant structures and slow dissociation. In Section \ref{subsec:res th0} we already noted that the homoclinic tangle of $\Gamma^i_+$ is responsible for a fractal structure of slow dissociation of initial conditions in $B_1^+$. Furthermore the homoclinic tangle of $\Gamma^o_+$ (and $\Gamma^o_-$) is responsible for slow dissociation in the interaction region $B_2$, especially at the top half of the energy interval.

  It is important to say that the section $\theta=const$, $p_\theta>0$ is not very well suited for the study of invariant manifolds. This is mainly due to the transition from vibration to rotation. The invariant manifolds $W_{\Gamma^i_+}$ may be nicely visible, but during this transition the invariant manifolds are not barriers to transport of surface area on this surface of section. Because parts $W_{\Gamma^i_+}$ rotate with $p_\theta<0$ after the transition, they do not return to the surface of section. For the same reason there are trajectories that do not return to the surface of section. The return map associated with this surface of section is therefore not area preserving. This anomaly can be seen from odd shapes of invariant manifolds - heteroclinic points seem to be mapped to infinity.

  Apart from the transition of $W_{\Gamma^i_+}$ from vibration to rotation, invariant manifolds may enter $B_1^\pm$ and be captured therein for a significant amount of time. Upon leaving $B_1^\pm$ the direction of rotation is unpredictable and this is true for invariant manifolds of all TSs. That is all we can say about the section $\theta=const$, $p_\theta>0$.

  The section on the inner DS, just as all other DSs, does not suffer from these problems, because they do not depend on the direction of rotation. Moreover, these surfaces are almost everywhere transversal to the flow.

\begin{figure}
  \centering
  {
  \includegraphics[width=5.6cm]{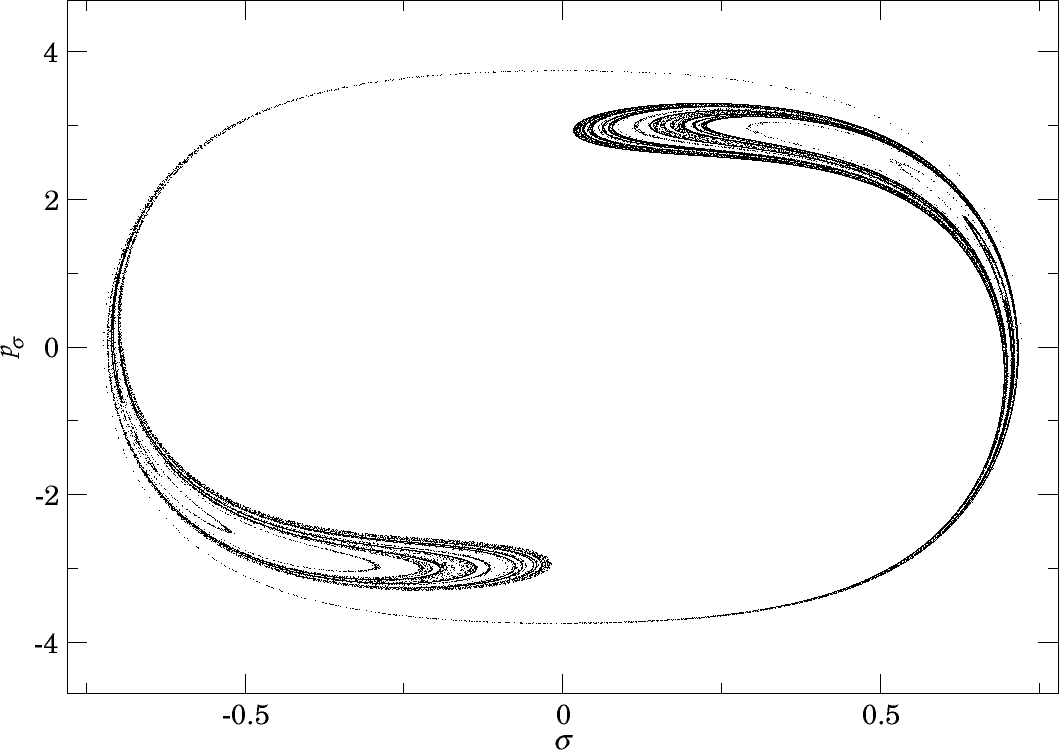}
  \includegraphics[width=5.6cm]{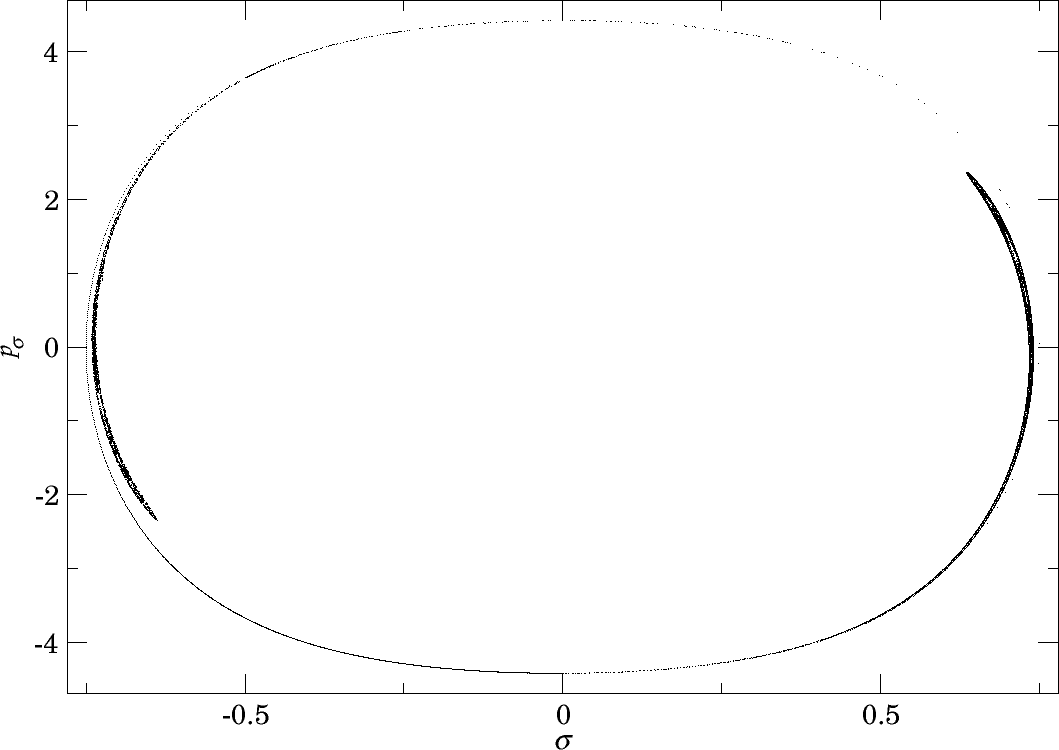}
  } 
  \caption{$W_{\Gamma^o_+}^{s-}$ invariant manifolds on the inner DS for $E=1$ (left) and $E=2$ (right). For energies $E\geq 2.5$ the manifolds $W_{\Gamma^o_+}^{-}$ don't reach the inner DS.}
  \label{fig:sec_disI}
  \end{figure}

  In Figure \ref{fig:sec_disI} we present the intersection of $W_{\Gamma^o_+}$ with the inner DS at $E=1$ and $E=2$. Since slow dissociation fades away at higher energies, we do not present the section at higher energies. In fact, for $E\geq 2.5$ the manifolds $W_{\Gamma^o_+}$ do not intersect the inner DS and therefore $W_{\Gamma^i_+}$ and $W_{\Gamma^o_+}$ do not intersect at all. Clearly then, slow dissociation, and thereby roaming, is induced by the heteroclinic tangle of $W_{\Gamma^i_+}$ and $W_{\Gamma^o_+}$. This claim is further supported by what we see in Figures \ref{fig:disI_rot} and \ref{fig:sec_disI}.

  When we compare Figures \ref{fig:disI_rot} and \ref{fig:sec_disI}, we clearly see that longer residence times are prevalent in the same locations where $W_{\Gamma^o_+}^{s-}$ intersects the inner DS. At $E=1$ we can even recognize the structure of the of the intersection in both figures. As the manifolds $W_{\Gamma^o_+}^{s-}$ recede with increasing energy, the area of slow dissociation at $E=2.5$ remains as a relic of the intersection. Afterall, trajectories close to $W_{\Gamma^o_+}^{s-}$ follow the manifold and approach $\Gamma^o_+$ before dissociation is completed.

  Note that there was no word of $W_{\Gamma^o_+}^{u-}$. This is mainly due to the fact that it influences the residence time in backward time, hence cannot be seen in forward time. Furthermore, $W_{\Gamma^o_+}^{u-}$ first intersects the other hemisphere of the inner DS, spends considerable time in $B_1^+$ and becomes heavily distorted before intersecting the outward hemisphere of the inner DS. In backward time, however, we expect a result symmetric to the one presented here due to time reversibility of the system.

  What really prevents us from making more fundamental conclusions at this point is the fact that $W_{\Gamma^o_+}^{s-}$ is heavily distorted when it reaches the inner DS. The reason is very simple - heteroclinic points. Here we not only mean trajectories on $W_{\Gamma^o_+}^{s-}$ that tend toward $\Gamma^i_+$, but also to $\Gamma^a_+$. Due to this fact, it is impossible to tell which area is enclosed by $W_{\Gamma^o_+}^{s-}$  and which is outside of it. For the majority of the area we note that $\sigma=0$, $p_\rho>0$, $p_\sigma=0$ (equivalent to $\theta=0$, $p_r>0$, $p_\theta=0$), the prototype of a fast dissociation, must lie inside $W_{\Gamma^o_+}^{s-}$ to quickly reach the outer DS. The tongues of $W_{\Gamma^o_+}^{s-}$ visible in Figure \ref{fig:sec_disI} therefore mostly contain trajectories that do not dissociate immediately.

  This problem is present on both the inner and outer DSs. Sections on both suffer from the fractal structure that is so characteristic for homoclinic and heteroclinic tangles. The ideal choice seems to be $\Gamma^a_+$ because $W_{\Gamma^i_+}^{+}$ and $W_{\Gamma^o_+}^{-}$ reach the middle DS very quickly. The first image of the manifolds under the Poincar\'e map associated with this surface does not display heteroclinic orbits, all manifolds are mapped to (topological) circles.

  Heteroclinic points become visible after applying the return map at least once, the resulting tongues wind around the previously mentioned circles. On the downside, $\Gamma^a_+$ is only hyperbolic until $2.72$, therefore for higher energies the middle DS allows local recrossings. It can be still used as a surface of section and we can expect to see fewer heteroclinic points that cause tongues, but we need to keep local recrossings in mind. In the next section we present a detailed view on the dynamics on the middle DS.


 \section{The observed dynamics and roaming}\label{sec:discussion}

  In this section we recall possible definitions of roaming used in previous works. We then elaborate on the observations above and analyse invariant manifolds on the middle DS with the aim to thoroughly explain how exactly roaming is linked to the heteroclinic tangles. Based on the explanation, a natural definition of roaming follows.

 \subsection{Roaming}
  Roaming in the chemistry literature refers to a kind of dissociation that is longer or more complicated than the usual dissociation with a monotonically increasing reaction coordinate that involves a saddle type equilibrium. While there is a sufficient amount of observations and intuitive understanding of what roaming is, an exact definition has not yet been generally adopted.

  Maugui\`ere et al. \cite{Mauguiere2014} proposed a classification of trajectories based on the number of turning points of trajectories in the interaction region $B_2$. Later the authors refine their definition in \cite{Mauguiere2014b} based on the number of intersections of a trajectory with the middle DS. Dissociating trajectories need to cross the middle DS at least three times before they are classified as roaming.

  Huston et al. \cite{Huston2016}, on the other hand, set the criteria such that roaming trajectories have to spend a certain amount of time at a minimum radius, have low average kinetic energy and have on average a certain number of bonds over time.


 \subsection{The mechanism of roaming}\label{subsec:roaming mech}

  Based on intersections of invariant manifolds, we would like to report on the types of trajectories in this dissociation problem and explain why the types exist. There is a general accord on the mechanism behind direct dissociation along the radical and molecular channel. The framework, that describes how codimension-$1$ invariant manifolds divide the energy surface in two and thereby separate reactive trajectories from non-reactive ones, is very well known in reaction dynamics, see \cite{MacKay84}, \cite{OzoriodeAlmeida90}, \cite{Rom-Kedar90}, \cite{Meiss15}. 

  Due to the different local geometries of the energy surface, we need to be careful with the invariant manifolds at this point. TSs that are brake orbits give rise to spherical DS and their invariant manifolds are spherical cylinders. TSs that are rotating orbits, just like ones belonging to the families $\Gamma^o$ and $\Gamma^a$, give rise a toric DS that is based on two orbits instead of one. Therefore in the description of transport, invariant manifolds of both orbits have to make up a toric cylinder together. Invariant manifolds govern transport of energy surface volume as follows.

  In $ \text{CH}_4^+ \rightarrow \text{CH}_3^+ + \text{H}$, we cannot discuss the molecular channel, but the radical channel and roaming is present. In general, if the H atom has enough kinetic energy to break bonds with $\text{CH}_3^+$, it escapes. Such a trajectory is contained in the interior of the invariant cylinder $W_{\Gamma^i_+}^{u+}$, because it leaves the inner DS to the $\text{CH}_3^++\text{H}$ side. The same is true for $W_{\Gamma^i_+}^{u+}$.

  Since the trajectory corresponding to $\theta=0$, $p_r>0$, $p_\theta=0$ on the inner DS dissociates immediately, a part of $W_{\Gamma^i_+}^{u+}$ reaches the middle and outer DS without returning to the inner DS. A part of the interior of $W_{\Gamma^i_+}^{u+}$ must therefore be contained in the invariant toric cylinder made up of $W_{\Gamma^a_+}^{s-}$ and $W_{\Gamma^a_-}^{s-}$, that we will refer to as $W_{\Gamma^a}^{s-}$. Other invariant toric cylinders will be denoted analogously.

  Trajectories that have too little energy in the radial degree of freedom do not reach the middle DS and are therefore not contained in the invariant cylinder. It does not matter whether $p_\theta>0$ or $p_\theta<0$. Considering that invariant manifolds are of codimension-$1$ on the energy surface and that $W_{\Gamma^a_+}^{s-}$ and $W_{\Gamma^a_-}^{s-}$ never intersect, by the inside of the invariant toric cylinder $W_{\Gamma^a}^{s-}$ we mean the energy surface volume enclosed between $W_{\Gamma^a_+}^{s-}$ and $W_{\Gamma^a_-}^{s-}$. As we shall see, $W_{\Gamma^i_+}^{u+}$ is entirely contained in the invariant cylinder $W_{\Gamma^a}^{s-}$.

  After crossing the middle DS, the interior of the cylinder $W_{\Gamma^a}^{s-}$ is lead away from the surface by the cylinder consisting of $W_{\Gamma^a}^{u+}$. The trajectories that dissociate are further guided by $W_{\Gamma^o}^{s-}$ towards the outer DS and further away by $W_{\Gamma^o}^{u+}$ to complete dissociation.

  All directly dissociating trajectories will be contained in the interior all of the above mentioned invariant cylinders. Moreover, directly dissociating trajectories are not contained in the interior of any other invariant cylinder.

  As soon as a trajectory is contained in another cylinder, it is guided by that cylinder to cross the corresponding DS. Should a trajectory be contained in $W_{\Gamma^i_+}^{s+}$, it will come back to the inner DS. In this way isomerisation, i.e. transport of energy surface volume between $B_1^+$ and $B_1^-$, is possible via the intersection of the interiors of $W_{\Gamma^i_+}^{u+}$ and $W_{\Gamma^i_-}^{s+}$ or $W_{\Gamma^i_-}^{u+}$ and $W_{\Gamma^i_+}^{s+}$.

  The intersections of the interiors of $W_{\Gamma^a}^{s+}$ and $W_{\Gamma^a}^{u+}$ or $W_{\Gamma^a}^{s-}$ and $W_{\Gamma^a}^{u-}$, on the other hand, lead to the recrossing of the middle DS. In case a trajectory originating in $B_1^\pm$ dissociates after recrossing of the middle DS, by the definition of Maugui\`ere et al. \cite{Mauguiere2014b} it is a roaming trajectory. From the above it is clear that roaming trajectories are contained in intersection of the interiors of $W_{\Gamma^i_\pm}^{u+}$, $W_{\Gamma^a}^{s-}$, $W_{\Gamma^a}^{u+}$, $W_{\Gamma^a}^{s+}$, $W_{\Gamma^a}^{u-}$ and $W_{\Gamma^o}^{s-}$. It remains to express the order of intersections of the DSs by a roaming trajectory with the invariant cylinders above.

  In summary the arguments above enable us to say that,
  \begin{itemize}
  
  \item directly dissociating trajectories are contained in $W_{\Gamma^i_+}^{u+}$ (or $W_{\Gamma^i_-}^{u+}$), $W_{\Gamma^a}^{s-}$, $W_{\Gamma^a}^{u+}$, $W_{\Gamma^o}^{s-}$ and no other,
  
  \item isomerisation and non-dissociating trajectories are contained in $W_{\Gamma^i_\pm}^{u+}$ and $W_{\Gamma^i_\mp}^{s+}$,
  
  \item roaming trajectories are contained in $W_{\Gamma^i_\pm}^{u+}$, $W_{\Gamma^a}^{s-}$, $W_{\Gamma^a}^{u+}$, $W_{\Gamma^a}^{s+}$, $W_{\Gamma^a}^{u-}$ and $W_{\Gamma^o}^{s-}$.
  
  \end{itemize}

  Note that since a trajectory contained in the cylinder $W_{\Gamma^a}^{s-}$ is automatically conveyed to $W_{\Gamma^a}^{u+}$ after crossing the middle DS, we may omit mentioning one of the cylinders. A roaming trajectory could therefore be shortly characterized by $W_{\Gamma^i_\pm}^{u+}$, $W_{\Gamma^a}^{s+}$ and $W_{\Gamma^o}^{s-}$.

  The definition of Maugui\`ere et al. admits nondissociating roaming trajectories. These are contained in $W_{\Gamma^i_\pm}^{u+}$ and $W_{\Gamma^a}^{s+}$, but not in $W_{\Gamma^o}^{s-}$.

  \begin{figure}
  \centering
  \includegraphics[width=10cm]{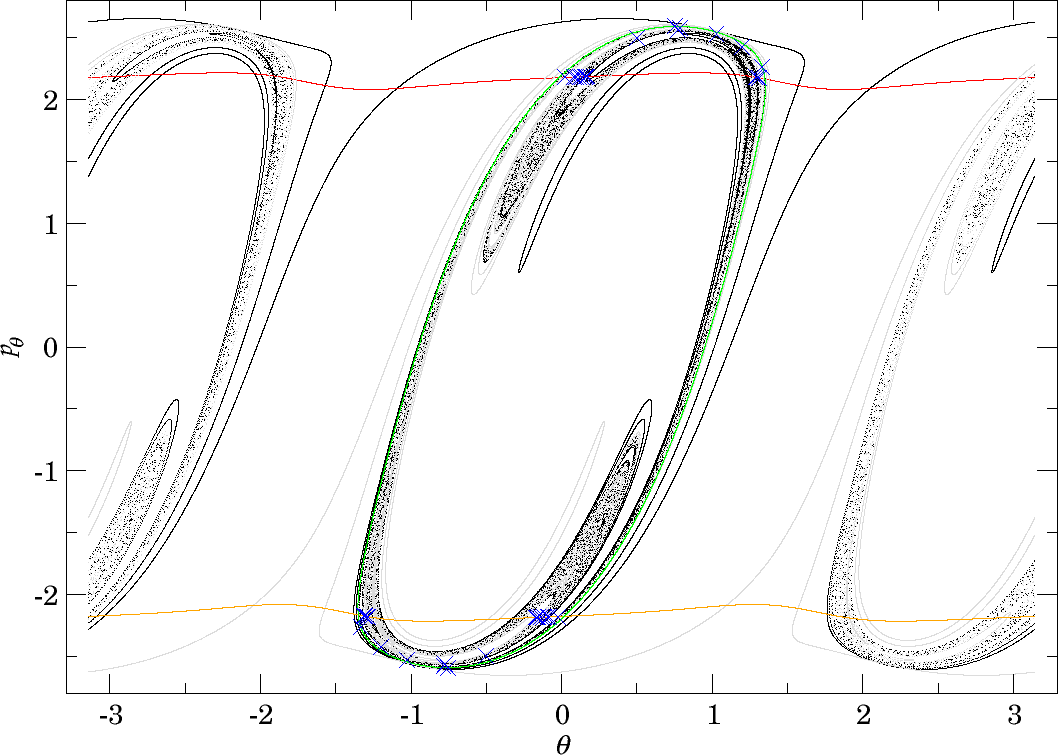}
  \caption{First and last intersections of invariant manifolds with the outward annulus of the middle DS for $E=1$. $W_{\Gamma^i_+}^{u+}$ (green) forms the boundary of $\gamma^{u+}_{i}$, $W_{\Gamma^o_+}^{s-}$ (red) and $W_{\Gamma^o_-}^{s-}$ (orange) form the boundary of $\gamma^{s-}_{o}$, $W_{\Gamma^o_+}^{u-}$ is black and $W_{\Gamma^o_-}^{u-}$ is grey. $W_{\Gamma^i_+}^{u-}$ copies the shape of $W_{\Gamma^o_+}^{u-}$ inside $\gamma^{u+}_{i}$. Selected initial conditions for roaming with very long residence times are marked with blue crosses.}
  \label{fig:roam_init1}
  \end{figure}

 \subsection{Roaming on the middle DS}\label{subsec:roaming middle DS}

  As mentioned in Section \ref{subsec:sec manifs}, the middle DS seems to be better suited for the study of roaming than the inner and outer DSs. More precisely, we will study dynamics on the outward annulus of the middle DS, i.e. the annulus crossed by the prototypical dissociating trajectory $\theta=0$, $p_r>0$, $p_\theta=0$. We may introduce canonical coordinates on this annulus using a generating function in the same way as we did in Section \ref{subsec:res inner DS}, but for for the sake of simplicity we continue using the coordinates $(\theta,p_\theta)$.

  In the following elaboration we need means to precisely express the order in which invariant cylinders intersect the outward annulus of the middle DS. Based on the arguments in Section \ref{subsec:roaming mech}, roaming involves the invariant cylinders $W_{\Gamma^i_\pm}^{u+}$, $W_{\Gamma^a}^{s+}$ and $W_{\Gamma^o}^{s-}$. Due to symmetry we have that every statement regarding $W_{\Gamma^i_+}^{u+}$ also holds for $W_{\Gamma^i_-}^{u+}$.

  The dynamics under the return map associated with the surface of section does not require $W_{\Gamma^a}^{s+}$ for a complete and detailed description of dynamics. The simple fact that a point on the surface is mapped by the return map to another point on the surface is enough to deduce that the corresponding trajectory is contained $W_{\Gamma^a}^{s+}$ and in fact, all other invariant cylinders made up of invariant manifolds of $\Gamma^a_\pm$.

  Consequently, for a description of roaming on the outward annulus of the middle DS we only need $W_{\Gamma^i_+}^{u+}$ and $W_{\Gamma^o}^{s-}$. Every branch of the invariant manifolds intersects the middle DS in a topological circle. Since it is possible that a branch of invariant manifold returns to the middle DS, by the first intersection of an unstable branch of invariant manifold with the outward annulus of the middle DS we mean that all points on the circle converge in backward time to the respective TS without reintersecting the outward annulus of the middle DS. Similarly we define the last intersection of a stable branch in forward time.

  \begin{figure}
  \centering
  \includegraphics[width=10cm]{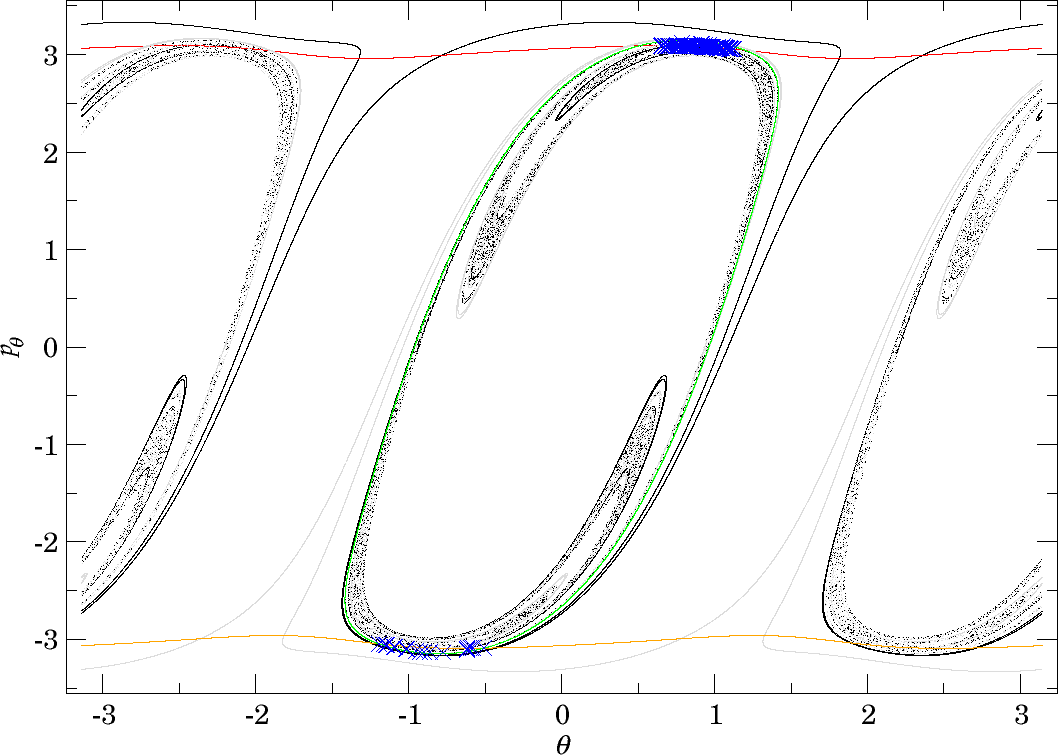}
  \caption{First and last intersections of invariant manifolds with the outward annulus of the middle DS for $E=2$. $W_{\Gamma^i_+}^{u+}$ (green) forms the boundary of $\gamma^{u+}_{i}$, $W_{\Gamma^o_+}^{s-}$ (red) and $W_{\Gamma^o_-}^{s-}$ (orange) form the boundary of $\gamma^{s-}_{o}$, $W_{\Gamma^o_+}^{u-}$ is black and $W_{\Gamma^o_-}^{u-}$ is grey. $W_{\Gamma^i_+}^{u-}$ copies the shape of $W_{\Gamma^o_+}^{u-}$ inside $\gamma^{u+}_{i}$. Selected initial conditions for roaming with very long residence times are marked with blue crosses.}
  \label{fig:roam_init2}
  \end{figure}

  Denote the interior of the first/last intersection of the invariant cylinders $W_{\Gamma^i_+}^{u+}$ and $W_{\Gamma^o}^{s-}$ with the outward annulus of the middle DS by $\gamma^{u+}_{i}$ and $\gamma^{s-}_{o}$, respectively. Denote the Poincar\'e return map associated with the outward annulus of the middle DS by $P$. By our findings all trajectories originating in $B_1^+$ and all trajectories that cross the inward annulus of the outer DS reach the middle DS.

\begin{figure}
  \centering
  \includegraphics[width=10cm]{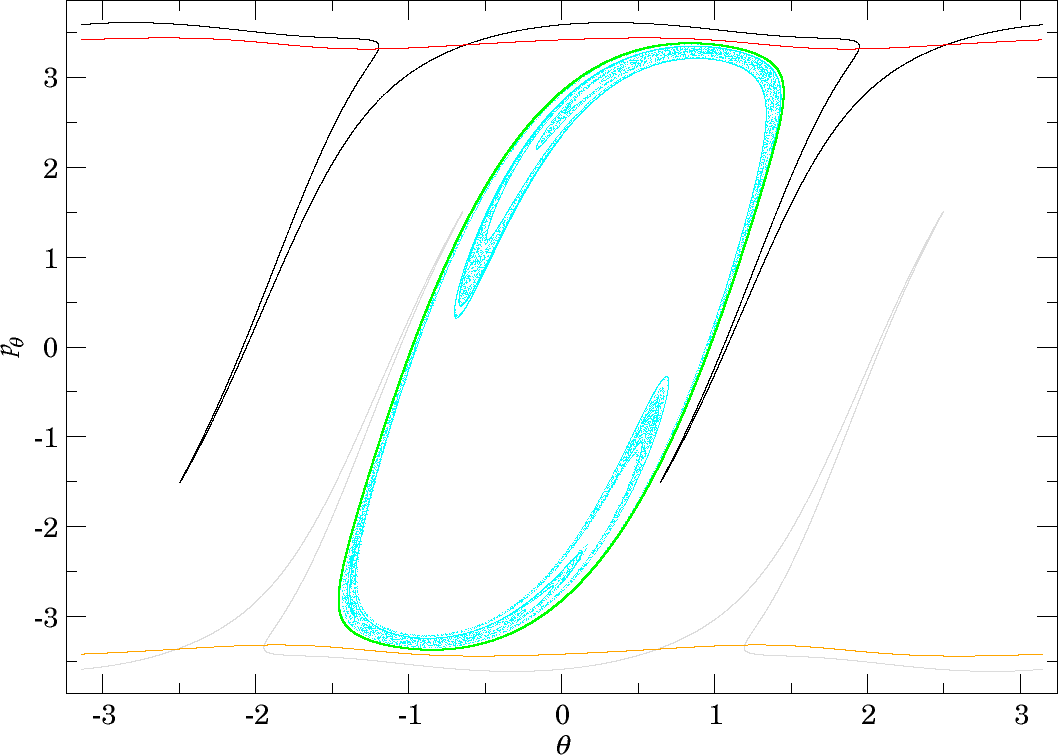}
  \caption{First and last intersections of invariant manifolds with the outward annulus of the middle DS for $E=2.5$. $W_{\Gamma^i_+}^{u+}$ (green) forms the boundary of $\gamma^{u+}_{i}$, $W_{\Gamma^o_+}^{s-}$ (red) and $W_{\Gamma^o_-}^{s-}$ (orange) form the boundary of $\gamma^{s-}_{o}$, $W_{\Gamma^o_+}^{u-}$ is black and $W_{\Gamma^o_-}^{u-}$ is grey, $W_{\Gamma^i_+}^{u-}$ is cyan. Roaming is not present because $\gamma^{u+}_{i}\subset \gamma^{s-}_{o}$ and $W_{\Gamma^i_+}^{u+}$ and $W_{\Gamma^o_+}^{s-}$ are disjoint.}
  \label{fig:roam_init2_5}
  \end{figure}

  By definition we have that $$\gamma^{u+}_{i}\cap\gamma^{s-}_{o},$$ contains trajectories that dissociate quickly. This is due to the fact that $\gamma^{u+}_{i}$ contains trajectories that just escaped from $B_1^+$ and $\gamma^{s-}_{o}$ contains those that reach the outer DS and therefore never return to the middle DS. Therefore points in $\gamma^{u+}_{i}\cap\gamma^{s-}_{o}$ do not have an image under the return map $P$, in fact the whole of $\gamma^{s-}_{o}$ does not have an image. This is in accordance with the results on ``reactive islands'' by \cite{OzoriodeAlmeida90}. Figures \ref{fig:roam_init1}, \ref{fig:roam_init2} and \ref{fig:roam_init2_5} show this intersection for various energies together with the first/last intersections of other invariant cylinders.

  Note that trajectories passing through $\gamma^{u+}_{i}\cap\gamma^{s-}_{o}$ reach the outer DS in varying amounts of time. We can expect the trajectory representing fast dissociation passing through $\theta=0$, $p_r>0$, $p_\theta=0$ to take significantly less time than trajectories in the proximity of $W_{\Gamma^o}^{s-}$, which may take arbitrarily long as they approach $\Gamma^o_\pm$.

  Therefore if roaming was to be only defined by time spent in $B_2$ or in the neighbourhood of a periodic orbit, we can always find a suitable trajectory in $\gamma^{u+}_{i}\cap\gamma^{s-}_{o}$ that is monotonous in $r$. Arguably, such a trajectory does not lead to an intramolecular hydrogen abstraction that has been reported in the context of roaming.

  It remains to explain what happens to $\gamma^{u+}_{i}\setminus\gamma^{s-}_{o}$. Since trajectories corresponding to points in $\gamma^{u+}_{i}\setminus\gamma^{s-}_{o}$ do not dissociate, they return to the outward annulus of the middle DS unless they are asymptotic to a periodic orbit. The set $\gamma^{u+}_{i}\setminus\gamma^{s-}_{o}$ has an image under the return map $P$ and it is $P\gamma^{u+}_{i}$. Note that the corresponding trajectories are guided to and from the middle DS by the invariant cylinders $W_{\Gamma^a}^{s-}$, $W_{\Gamma^a}^{u+}$, $W_{\Gamma^a}^{s+}$, $W_{\Gamma^a}^{u-}$. We remark that based on the understanding of lobe dynamics (\cite{Rom-Kedar90}, trajectories that cross the outward annulus of the middle DS and do not reach the outer DS are repelled by $\Gamma^o_\pm$ and necessarily pass through the homoclinic tangle of $\Gamma^o_\pm$.

  By the definition of Maugui\`ere et al. \cite{Mauguiere2014b}, roaming trajectories cross the middle DS at least three times, which means crossing the outward annulus at least twice and the inward annulus at least once. Roaming trajectories must therefore contained in $P\gamma^{u+}_{i}$. In fact roaming trajectories that cross the outward annulus of the middle DS precisely $n$ times before dissociating pass through
  $$P^{n-1}\gamma^{u+}_{i}\cap \gamma^{s-}_{o}.$$
  Since recrossings of the middle DS are possible due to the homoclinic tangle of $\Gamma^o_\pm$, roaming requires that the invariant cylinder $W_{\Gamma^i_+}^{u+}$ conveys trajectories into the homoclinic tangle. Heteroclinic intersections are therefore necessary.

  Arguably, recrossings of the middle DS are inevitable to capture the process of intramolecular hydrogen abstraction, as reported by \cite{Bowman2011}, where the free H atom has to return back to the $\text{CH}_3^+$ core.

  Isomerisation trajectories are also contained in $\gamma^{u+}_{i}\setminus\gamma^{s-}_{o}$ and are guided by $W_{\Gamma^i_\pm}^{s+}$ to $B_1^\pm$ and by $W_{\Gamma^i_\pm}^{u+}$ out of $B_1^\pm$. Trajectories that return to $B_1^+$ pass through the intersection of $P^n\gamma^{u+}_{i}\cap W_{\Gamma^i_-}^{u+}$, for some $n$. Trajectories corresponding isomerisation pass through the intersection of $P^n\gamma^{u+}_{i}$ and the last intersection of $W_{\Gamma^i_-}^{s+}$ with the outward annulus, for some $n$.


  It remains to discuss the first intersection of $W_{\Gamma^o_+}^{u-}$, $W_{\Gamma^o_-}^{u-}$ and $W_{\Gamma^i_+}^{u-}$ on the surface of section shown in Figures \ref{fig:roam_init1}, \ref{fig:roam_init2} and \ref{fig:roam_init2_5}. We shall denote the intersections according to the convention above by $\gamma^{u-}_{o}$ and $\gamma^{u-}_{i}$, respectively.

  The invariant cylinder $W_{\Gamma^o}^{u-}$ guides energy surface volume from the inward annulus of the outer DS into the interaction region. Clearly then $\gamma^{s-}_{o}\cap \gamma^{u-}_{o}$ corresponds to trajectories that intersect the surface of section only once. The part of the intersection in $\gamma^{u+}_{i}$ passes through $B_1^+$. $W_{\Gamma^o}^{u-}$ is guided from the inward hemisphere of the inner DS by $W_{\Gamma^i_+}^{u-}$ and its homoclinic intersections cause tongues. In the process, $W_{\Gamma^o}^{u-}$ and $W_{\Gamma^i_+}^{u-}$ are stretched and compressed causing only one to be visible in Figures \ref{fig:roam_init1} and \ref{fig:roam_init2}. In Figure \ref{fig:roam_init2_5} $W_{\Gamma^o}^{u-}$ does not enter $B_1^+$ and the two invariant cylinders are visible.

  It is important to point out that seemingly $W_{\Gamma^o}^{u-}$ and $W_{\Gamma^i_+}^{u-}$ intersect, which is impossible. Instead we observe a discontinuity caused by points on $W_{\Gamma^o}^{u-}$ heteroclinic to $\Gamma^i_+$. As mentioned above, a part of $W_{\Gamma^o}^{u-}$ passes through $B_1^+$ and is mapped by $P$ into $\gamma^{u+}_{i}$, the remainder visible on the surface of section stays in $B_2$ and is mapped outside of $\gamma^{u+}_{i}$. The points inbetween are not intersections between $W_{\Gamma^o}^{u-}$ and $W_{\Gamma^i_+}^{u-}$, but between $W_{\Gamma^o}^{u-}$ and $W_{\Gamma^i_+}^{s\pm}$ and do not have an image under $P$.

  Note that $\gamma^{u-}_{o}$ carries information about roaming in backward time. Since points in $\gamma^{u-}_{o}$ do not have a preimage under $P$,
  all points in $\gamma^{u+}_{i}\cap\gamma^{s-}_{o}$ that are not in $\gamma^{u-}_{o}$ must have a preimage under $P$. These points correspond to trajectories that qualify as roaming in backward time. With increasing energy it becomes difficult to study $\gamma^{u+}_{i}\setminus\gamma^{s-}_{o}$ due to the fractal structures of $W_{\Gamma^o}^{u-}$ caused by heteroclinic points. We can, however, expect proportionally fewer roaming trajectories to enter $B_1^+$ multiple times at $E=2$ than at $E=1$. Instead it is probable to find roaming trajectories spending the majority of their residence time in $B_2$ at $E=2$.


 \section{Global study of the invariant manifolds that govern the dynamics}\label{sec:representation}

  In this section we discuss an alternative way of studying dynamics on a $3$-dimensional energy surface using the so called Conley-McGehee representation \cite{Conley1968}, \cite{McGehee1969}, \cite{MacKay1990}, described along with other alternatives in \cite{Waalkens2010}. This is a very useful way of studying dynamics in full $3$ dimensions, but to date has only been defined for subsets of energy surfaces that are locally a spherical shell. Since the Conley-McGehee representation does only works in $B_1^\pm$ of Chesnavich's CH$_4^+$ model studied here, we introduce an extension of the Conley-McGehee representation that enables us to study energy surfaces with other geometry than in the Conley-McGehee case.

 \subsection{Conley-McGehee representation}
  The dynamics on the energy surface can be visualized in many ways. Just as it was done above, it can be viewed on various surfaces of section, most notably ones constructed around TSs. It is also possible to study the system locally using normal form approximations. The Williamson normal form \cite{Williamson36} of the Hamiltonian in the neighbourhood of an index-$1$ critical point is
  $$H_2(q_1,p_1,q_2,p_2)=\frac{1}{2}\lambda (p_1^2-q_1^2)+\frac{1}{2}\omega (p_2^2+q_2^2),$$
  for some $\lambda, \omega>0$.
  We found that for a fixed energy $H_2(q_1,p_1,q_2,p_2)=h_2$, the energy surface can be locally viewed as a continuum of spheres parametrized by $q_1$.

  In the Conley-McGehee representation \cite{Conley1968}, \cite{McGehee1969}, is based on the spherical local geometry of an energy surface. While the normal form perspective above only applies locally, in the Conley-McGehee representation the whole energy surface is represented as a nested set of spheres parametrised in the radial direction by the reaction coordinate.

  Advantages are immediate - the representation gives a global model of the energy surface and by construction reveals the spherical structure of the energy surface. For $2$ degrees of freedom it enables us to study the $3$-dimensional energy surface in the full $3$ dimensions. Moreover, it enables to visualise the DSs as spheres that separate the energy surface into two disjoint components. It is also very natural that the flux through the hemispheres of the DSs is unidirectional and trajectories have to cross a particular hemisphere of the DS to pass from one component to the other.

  Apart from DSs the Conley-McGehee representation enables us to visualise and therefore study TSs and their invariant manifolds that are spherical cylinders in a natural environment.

 \subsection{Toric extension of the Conley-McGehee representation}
  The Conley-McGehee representation in its original form applies to spherical geometries. The energy surface of the CH$_4^+$ model, on the other hand, has a partially spherical and partially toric geometry, where many periodic orbits come in pairs and several DSs are tori. We therefore adapt the Conley-McGehee representation for the energy surface as follows.

  The energy surface is defined by
  $$ M_E = \Big\lbrace (r,\theta,p_r,p_\theta) \in \Rbb^4 \Big\vert H(r,\theta,p_r,p_\theta)=\frac{1}{2 m} p_r^2 + \frac{1}{2} \left(\frac{1}{I} + \frac{1}{m r^2} \right) p_\theta^2 + U(r,\theta) = E \Big\rbrace.$$
  For very high energies, $E>E_2$, where only $r<0.9$ is energetically inaccessible due to the cut-off of the potential and $E > U(r,\theta)$ for all $r\geq0.9$, the whole energy surface has a toric local geometry. For any fixed radius $r_0$ and a fixed $\theta_0$,
  $$M_E(r_0,\theta_0)= \Big\lbrace (p_r,p_\theta) \in \Rbb^2 \Big\vert \frac{1}{2 m} p_r^2 + \frac{1}{2} \left(\frac{1}{I} + \frac{1}{m r_0^2} \right) p_\theta^2 = E-U(r_0,\theta_0)\Big\rbrace,$$
  is a $\Sbb^1$. If $\theta$ is not fixed, $M_E(r_0,\theta)$ defines a $\Sbb^1\times\Sbb^1=\Tbb^2$ and therefore the whole energy surface $M_E$ is a $\Tbb^2\times \Rbb^+$.

  The radii of the concentric circles $M_E(r_0,\theta_0)$ on the $(p_r,p_\theta)$-plane depend on $r_0$ and $\theta_0$ through $U(r_0,\theta_0)$. The potential energy is not monotonous in $r$ nor in $\theta$. Recall from Section \ref{subsec:hill} that
  \begin{itemize}
  \item the two wells $q_0^\pm$ are located at $(1.1, 0)$ and $(1.1, \pi)$ with $U(q_0^\pm)=E_0\approx -47$,
  \item the two index-$2$ saddles $q_2^\pm$ are located at $(1.63, \frac{\pi}{2})$ and $(1.63, \frac{3\pi}{2})$ with $U(q_2^\pm)=E_2\approx 22.27$.
  \end{itemize}
  Further recall from Section \ref{subsec:general setting} that for $r$ sufficiently large $U(r,\theta)$ is essentially independent of $\theta$ and $r^{2}U(r,\theta)\rightarrow0$ as $r\rightarrow\infty$ for all $\theta$.

  It follows that the tori corresponding to $r_0=1.1$ and $r_0=r_{large}$, for some $r_{large}$ sufficiently large, intersect. This is because the circle $M_E(1.1,0)$ has a larger radius than $M_E(r_{large},0)$, while $M_E(1.1,\frac{\pi}{2})$ has a smaller radius than $M_E(r_{large},\frac{\pi}{2})$.

  The tori will always intersect if the radius is not a monotonous in $r$. In order to extend the Conley-McGehee representation, we need to reparametrise these tori so that their radii are monotonous in $r$ for every $\theta$.

  Define
  $$ P_r = \frac{r}{\sqrt{2 m (E - U(r,\theta)) } } p_r, $$
  and
  $$ P_\theta = \frac{r \sqrt{\frac{1}{I} + \frac{1}{m r^2} }}{\sqrt{2  (E - U(r,\theta)) } } p_\theta.$$
  Now we have
  $$ P_r^2+P_\theta^2 = r^2\frac{1}{E - U(r,\theta)} \Big(\frac{1}{2 m} p_r^2+\frac{1}{2}\left(\frac{1}{I} + \frac{1}{m r^2}\right) p_\theta^2 \Big)=r^2.$$

  The radius of the tori is monotonous in $r$ and independent of $\theta$ and therefore the tori $P_r^2+P_\theta^2=r^2$ foliating the energy surface are, unlike before the reparametrisation, disjoint in $(\theta,P_r,P_\theta)$-space.

  Note that a section of the tori with a plane of section $\theta=\theta_0$ shows concentric circles, where the smallest one has the radius $r=0.9$ due to the cut-off of the system. This is due to the fact that the boundary of the energy surface corresponds to a torus. Should it be desirable to have the whole $(\theta,P_r,P_\theta)$-space foliated by tori, it can be done by replacing $r$ by $r-0.9$ in the definitions of $P_r$ and $P_\theta$.

 \subsection{Extension to non-constant geometries}
  We remark that the construction above relies on the fact that at $E>E_2$ the energy surface has a purely toric geometry. We can slightly amend the construction to work for lower energies $E\leq E_2$, where the energy surface is not purely toric.

  $E = U(r,\theta)$ does not pose a problem for the definition of $P_r$ and $P_\theta$ as it may seem on first sight. $P_r$ and $P_\theta$ are only normalized conjugate momenta and by definition
  $$|P_r|,|P_\theta|\leq r.$$
  The momenta are therefore well defined on the whole energy surface.

  Points on $E = U(r,\theta)$ are degenerate circles with radius $0$ on the energy surface, but due to normalization correspond to circles in $(P_r,P_\theta)$. Such a representation of the energy surface for lower energies is clearly flawed. For $r$ large, we still have tori, but for smaller $r$, e.g. near $\Gamma^i$, we do not see the spherical geometry we expect.

  To solve this issue, we introduce different momenta in which the radius $P_r^2+P_\theta^2\rightarrow 0$ as $U(r,\theta)\rightarrow E$. We remark that the following only works for $E\leq \widetilde{E}_1$. In the interval $\widetilde{E}_1<E<E_2$, the projection of the energy surface on configuration space is the whole plane minus three discs, see Figure \ref{fig:potentialcontour1}. One is the potential energy cut-off and the other two are areas of high potential around index-$2$ critical points $q_2^\pm$. Spherical and toric geometry cannot accurately represent a genus $3$ surface.

  Since for energies $E<0$ the standard Conley-McGehee representation applies, we will restrict ourselves to the more interesting case $0\leq E\leq \widetilde{E}_1$.
  For the sake of simplicity, we retain the notation $P_r$ and $P_\theta$, making clear that we are discussing different momenta in a different energy interval than before.

  \begin{figure}
  \centering
  \includegraphics[width=.5\textwidth]{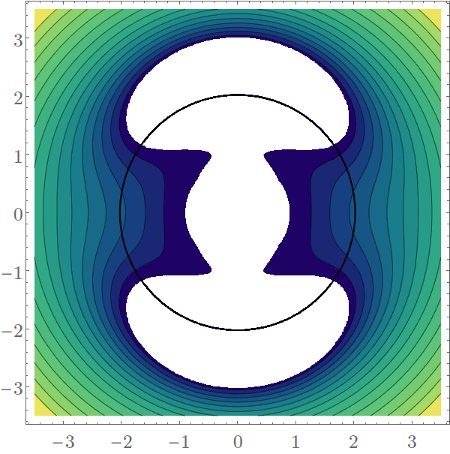}
  
  \caption{Contour plot of $r^6 \left( E - U(r,\theta) \right)$ for $E=0$. The black circle corresponds to the radius $r_E$ defined in the text.}
  \label{fig:ConleyContour}
  \end{figure}

  Let
  $$ P_r = \frac{r^3}{\sqrt{2 m} } p_r, $$
  and
  $$ P_\theta = r^3 \sqrt{\frac{1}{2}\left(\frac{1}{I} + \frac{1}{m r^2}\right) } p_\theta.$$
  It follows that
  $$ P_r^2+P_\theta^2 = r^6 \Big(\frac{1}{2 m} p_r^2+\frac{1}{2}\left(\frac{1}{I} + \frac{1}{m r^2}\right) p_\theta^2 \Big)=r^6 (E - U(r,\theta)).$$

  While $E - U(r,\theta)$ makes sure that zero kinetic energy corresponds to $P_r=P_\theta=0$, the term $r^6$ seems perhaps less obvious. In the previous section we showed that it is important for $P_r^2+P_\theta^2$ to be monotonous in $r$ for every $\theta$. This is however not possible, because for some fixed angles $\theta=const$ the term $(E - U(r,\theta))$ vanishes for several values of $r$. This is only possible in coordinates in which $E = U(r,\theta)$ correspond to coordinate lines for all $E$.

  Let $r_E$ be the smallest $r$ such that $E = U(r,\theta)$ has at most one solution for every $\theta$ for $r_E\leq r$. The term $r^6$ is the smallest even power of $r$ such that $r^6 (E - U(r,\theta))$ is monotonous in $r$ on $r_E\leq r$ for $0\leq E\leq \widetilde{E}_1$.

  As a consequence of restricting the radius, the representation omits a significant part of $B_1^\pm$. Since we formulated roaming as a transport problem from the inner DS to the outer DS, dynamics inside $B_1^\pm$ does not play a significant role in our study. All significant periodic orbits and DSs are well defined in the Conley-McGehee representation as presented here.

  Figure \ref{fig:ConleyContour} shows the contour plot of $r^6 (E - U(r,\theta))$ with a highlighted circle marking $r_0$, the boundary of the representation defined above.

 \subsection{Consequences of the extension}

  \begin{figure}
  \centering
  \includegraphics[width=.49\textwidth]{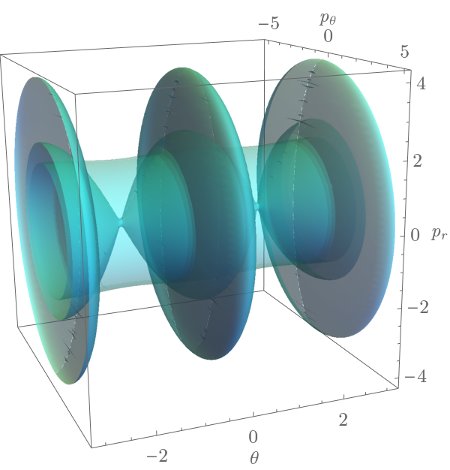}
  \includegraphics[width=.49\textwidth]{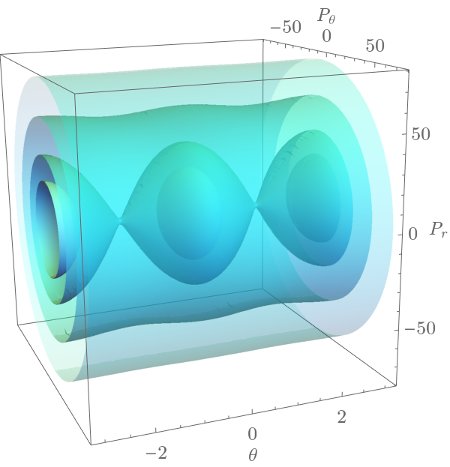}
  
  \caption{Comparison of energy surface geometry representation for $E=0$ in $(\theta,p_r,p_\theta)$, left, and $(\theta,P_r,P_\theta)$, right. The surfaces shown correspond to selected fixed values of $r$. The value $r=r_0=2.0267$ corresponds to the spheres, $r=2.802$ defines the pinched torus and $r=3.5$ is a regular torus. Additionally, the right figure also shows $r=4$.}
  \label{fig:representation}
  \end{figure}

  In the representation as defined above, the geometry of the energy surface is preserved. For a fixed radius $r_0$, we can see that the surfaces have the following topologies:
  \begin{itemize}
  \item $\Tbb^2$ if $U(r_0,\theta) < E$ for all $\theta$,
  \item a pinched torus if $U(r_0,\theta) \leq E$ for all $\theta$ and $U(r_0,\theta_0) = E$ for some $\theta_0$,
  \item $\Sbb^2 \cup \Sbb^2$ if $U(r_0,\theta_0) > E$ for some $\theta_0$.
  \end{itemize}

  For $E=0$, an example of each is shown in Figure \ref{fig:representation} in the canonical phase space coordinates $(\theta,p_r,p_\theta)$ and in the proposed extension of Conley-McGehee representation $(\theta,P_r,P_\theta)$. Note indeed in the latter that the surfaces are disjoint and present part of a foliation of the energy surface. We added an additional value of $r=4$ to the extended Conley-McGehee representation to illustrate that the radius of the tori diverges for $r\rightarrow\infty$, whereas it converges in the canonical phase space coordinates.

  Due to the properties of the extended Conley-McGehee representation, we may study invariant structures and the aforementioned DSs globally on the energy surface. All techniques used to date either relied on surfaces of section or local approximations of the energy surface. In what follows, we study the structures on the energy surface in full three dimensions.

\begin{figure}
  \centering
  \includegraphics[width=10cm]{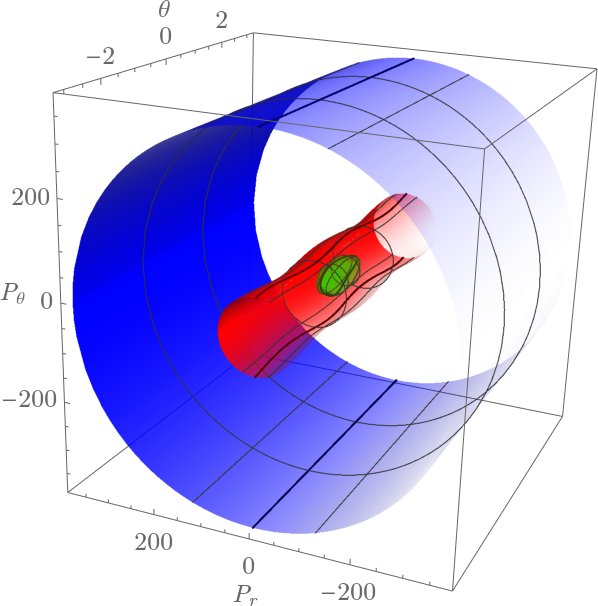}
  \caption{TSs and DSs in the Conley-McGehee representation for $E=2.5$. The inner DS is shown in green, middle DS in red, outer DS in blue. $\Gamma^i_+$, $\Gamma^a_\pm$ and $\Gamma^o_\pm$ are the thick lines on the corresponding DSs shown in darker green, red and blue, respectively.}
  \label{fig:ConTS25}
  \end{figure}

  In Figure \ref{fig:ConTS25} we present the periodic orbits $\Gamma^i_+$, $\Gamma^a_\pm$ and $\Gamma^o_\pm$, all of which are TSs, and the associated DSs at energy $E=2.5$. The construction of the DSs was discussed in Section \ref{subsec:DS}. We remark that from a qualitative perspective Figure \ref{fig:ConTS25} can be thought to represent the whole energy interval $0<E< 2.72$, where $\Gamma^a_\pm$ are unstable. The difference at higher energies is that $\Gamma^a_\pm$ is not a TS and the associated torus is not a DS.

  We have used some of the DSs mentioned above in previous sections to study residence times, rotation numbers, and most importantly, the intersections of stable and unstable invariant manifolds of TSs. For the sake of clarity, in the following we left out the $\Gamma^i_-$ and all associated structures, but everything said about $\Gamma^i_+$ also holds for $\Gamma^i_-$. In the figures one can easily imagine another sphere just like the inner DS but shifted by $\pi$ in the angular direction.

  Note that here we take full advantage of the proposed extension of the Conley-McGehee representation to present structures on the energy surface with different local geometries - the inner DS is a sphere whereas the middle and the outer DSs are tori. This has to our knowledge not been done before.

  On the DSs we highlighted the respective TSs. Note that the inner DS is defined using one periodic orbit whereas the middle and the outer are defined using two. The individual orbits in the families $\Gamma^a$ and $\Gamma^o$ can be distinguished by the sign of $P_\theta$ as they run in opposite directions.

  The inner DS is divided by $\Gamma^i_+$ into two hemispheres and the surfaces of unidirectional flux from the definition of a TS. Flux from $B_1^+$ to $B_2$ crosses the hemisphere where predominantly $P_r>0$. The situation is similar for the annuli of the middle and outer DS. The two orbits divide the torus into two annuli, where the outward annulus is the one with larger $P_r$.

  \begin{figure}
  \centering
  \includegraphics[width=10cm]{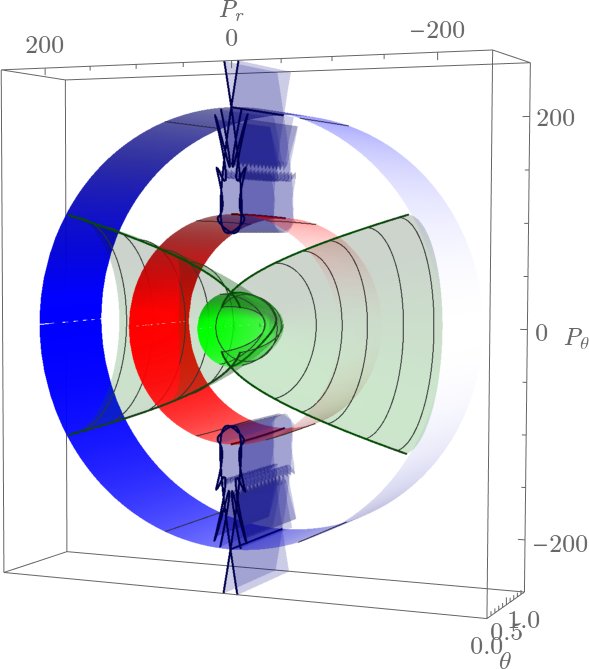}
  
  \caption{TSs, corresponding DSs and invariant manifolds in the Conley-McGehee representation at $E=5$. The inner DS is shown in green, middle DS in red and outer DS in blue. $\Gamma^i_+$, $\Gamma^a_\pm$ and $\Gamma^o_\pm$ are the thick lines on the corresponding DSs shown in darker green, red and blue respectively. Note that $\Gamma^a_\pm$ are stable at this energy. The invariant manifolds $W_{\Gamma^i_+}$ (green) are only a sketch based on the computed sections on the surface $\theta=0$ shown in thick green.}
  \label{fig:ConManifFull}
  \end{figure}

  The mechanism behind energy surface volume transport across DSs is governed by invariant manifolds of the corresponding TSs as discussed in the Sections \ref{subsec:sec manifs} and \ref{sec:discussion}. Here we present a new perspective for the study of invariant manifolds. In the following we use $E=5$, because at this energy the TSs are evenly spaced and the lack of heteroclinic intersections facilitates understanding of this new perspective. Everything we say is applicable to the whole interval $0<E<6.13$ relevant to roaming.

  Figure \ref{fig:ConManifFull} displays the invariant manifolds of $\Gamma^i_+$ and $\Gamma^o_\pm$ in full $3$ dimensions. Note that the manifolds $W_{\Gamma^i_+}$ are only a sketch based on computed sections on the surface $\theta=0$. Computing the whole invariant manifold numerically is in this case relatively straight-forward, but for the sole purpose of illustration unnecessarily expensive.

  \begin{figure}
  \centering
  \includegraphics[width=10cm]{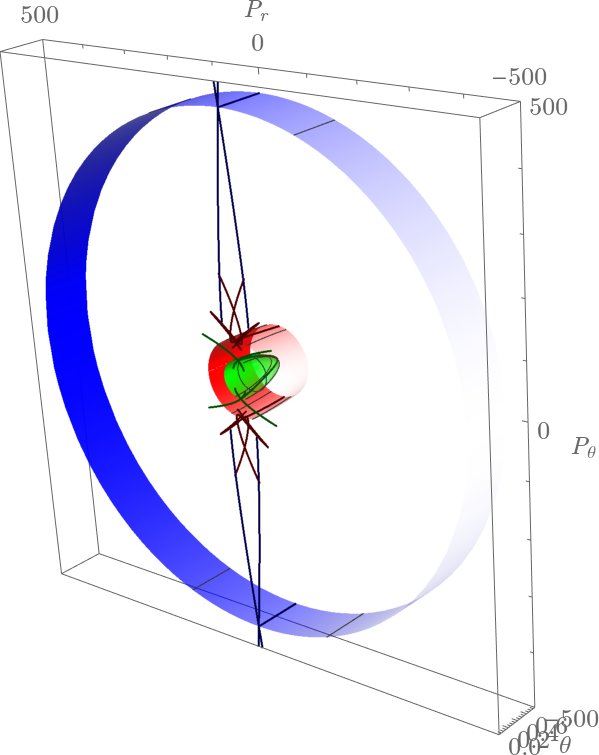}
  
  \caption{A section of invariant manifolds for $\theta=0$ at $E=2$ in the Conley-McGehee representation. The inner DS is shown in green, middle DS in red, outer DS in blue. $\Gamma^i_+$, $\Gamma^a_\pm$ and $\Gamma^o_\pm$ are the thick lines on the corresponding DSs shown in darker green, red and blue respectively and the curves extending from the TSs are their respective invariant manifolds. For clarity, the intersections of $W_{\Gamma^i}$ and $W_{\Gamma^o}$ are only indicated.}
  \label{fig:ConManif2}
  \end{figure}

\begin{figure}
  \centering
  \includegraphics[width=10cm]{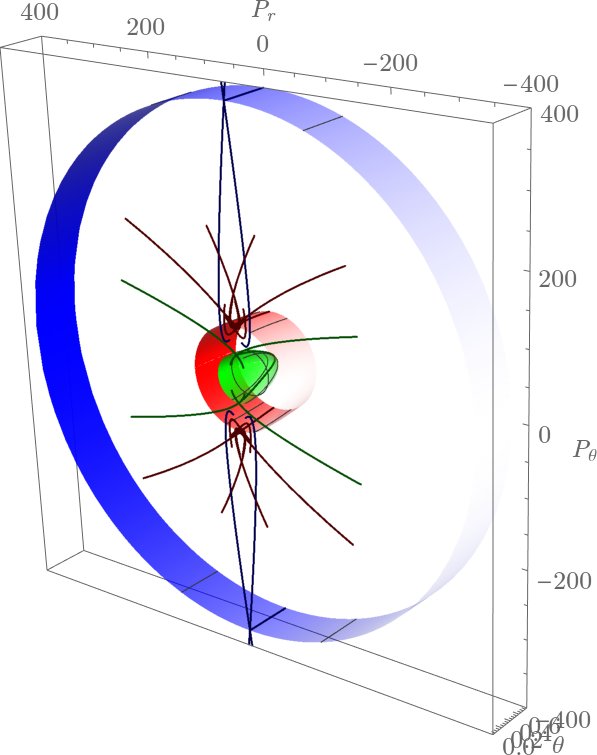}
  \caption{A section of invariant manifolds for $\theta=0$ at $E=2.5$ in the Conley-McGehee representation. The inner DS is shown in green, middle DS in red, outer DS in blue. $\Gamma^i_+$, $\Gamma^a_\pm$ and $\Gamma^o_\pm$ are the thick lines on the corresponding DSs shown in darker green, red and blue respectively and the curves extending from the TSs are their respective invariant manifolds. Note that $W_{\Gamma^i}$ and $W_{\Gamma^o}$ do not intersect.}
  \label{fig:ConManif25}
  \end{figure}

  Clearly visible is the structure of the manifolds, spherical cylinders formed by $W_{\Gamma^i_+}$ and toric cylinders formed by $W_{\Gamma^o_\pm}$. This is how the sections of invariant manifolds in the extended Conley-McGehee representation near $\theta=0$ shown in Figures \ref{fig:ConManif2} and \ref{fig:ConManif25} should be interpreted. For clarity, the invariant cylinders are indicated by the section of invariant manifolds for $\theta=0$ in these figures, the scale and complexity of intersections of the invariant cylinders would make the figures incomprehensible.

  One can clearly see that at $E=5$, in fact in the whole energy interval $E\geq2.5$, the intersection of $W_{\Gamma^i_+}^{u+}$ with the middle and the outer DSs produces a topological circle centred at $P_\theta=0$. This is purely the consequence of the spherical geometry induced by $\Gamma^i_+$. $W_{\Gamma^o_\pm}$ for the same reason intersects the middle DS in two lines that should be seen as circles concentric with $\Gamma^a_\pm$. This is in agreement with the sections of invariant manifolds on the outward annulus of the middle DS at $E=2.5$ shown in Figure \ref{fig:roam_init2_5}.
  
  Note that $W_{\Gamma^i_+}^{u-}$ intersects the outward hemisphere of the middle DS at $E=2.5$ in a shape that cannot be identified as a circle. This is due to the fact that we study $W_{\Gamma^i_+}^{u-}$ that is asymptotic to $\Gamma^i$ in $B_1^+$ as it leaves from $B_1^+$ to $B_2$. Moreover, it is deformed in the proximity of $\Gamma^a$ that exhibits a different kind of dynamics than $\Gamma^i$, rotating as opposed to vibrating.

  According to the findings in Section \ref{subsec:roaming middle DS}, the non-existence of roaming at higher energies due to the lack of intersection of $W_{\Gamma^i_\pm}$ with $W_{\Gamma^o_\pm}$ is immediate from Figures \ref{fig:ConManifFull} and \ref{fig:ConManif25}.

  The situation at $E<2.5$ represented by Figure \ref{fig:roam_init2} is very similar to the energy interval $E\geq2.5$. The main difference here are the heteroclinic intersections that cause roaming. These intersections are visible on the middle DS as well as in the extended Conley-McGehee representation. We remark that in the extended Conley-McGehee representation roaming occurs in the thin stripes between the two invariant cylinders $W_{\Gamma^o}^{s-}$ and $W_{\Gamma^o}^{u-}$ around $P_r=0$. Clearly the majority of the energy surface is occupied by more direct dynamics.


 \section{Conclusion}\label{sec:conclusions}
  We have shown that numerical observations of long dissociation are caused by particular structures formed by invariant manifolds of TSs. These invariant manifolds are also responsible for multiple recrossings of the middle DS and consequently also for roaming.

  We have shown that roaming trajectories that originate in the potential wells are captured in the homoclinic tangles of the outer TS and the middle TS dissociating. The transition for the potential wells into the homoclinic tangles is only possible in case the invariant manifolds of the inner and outer DS intersect and create a heteroclinic tangle. In case of Chesnavich's CH$_4^+$ model, this heteroclinic intersection only exists for energies $E\leq 2.5$ and therefore the system does not admit roaming at higher energy levels.

  Our results can possible be directly extended to other chemical reactions, as the only significant assumption on the potential energy is that for all $\theta$ $$U(r,\theta)\in o(r^{-2}) \text{ as } r\rightarrow\infty.$$
  This condition guarantees the existence of an outer TS thanks to which we may restrict roaming to a transport problem from potential wells representing a stable molecule to the DS associated with the outer TS.

  Furthermore, our findings support then dynamical definition of roaming as suggested in \cite{Mauguiere2014b}.

  We studied the above-mentioned invariant manifolds on various dividing surfaces, some of which highlighted the difficulties posed by the unusual energy surface geometry. This was the case even hough all surfaces of section satisfy the Birkhoff condition \cite{Birkhoff27} of being bounded by invariant manifolds. These difficulties motivated us to construct an extension of the Conley-McGehee representation. Using the extension we were able to study the unusual energy surface geometry, TSs, the associated invariant manifolds and DSs in full three dimensions.


 \section*{Acknowledgment}

  We are grateful to Dayal Strub for many very useful discussions and remarks on the contents of this paper.


\bibliographystyle{plain}
\bibliography{refs}

\end{document}